\documentclass[11pt,a4paper]{article}
\usepackage{amsmath}
\usepackage{amsfonts}
\usepackage{amssymb}
\usepackage{algorithm}
\usepackage{algorithmic}
\usepackage{booktabs}
\usepackage{color}
\usepackage{enumerate}
\usepackage{graphicx}
\usepackage[hidelinks]{hyperref}
\usepackage{lscape}
\usepackage{longtable}
\usepackage{mathrsfs}
\usepackage{rotating}
\usepackage{subfigure}
\usepackage{url}
\usepackage[dvipsnames]{xcolor}
\usepackage{booktabs}

\usepackage{titlesec} 
\usepackage[titletoc]{appendix}

\newtheorem{theorem}{Theorem}[section]

\newtheorem{corollary}{Corollary}[section]

\newtheorem{lemma}{Lemma}[section]

\newtheorem{remark}{Remark}[section]

\newtheorem{assumption}{Assumption}[section]

\newenvironment{proof}[1][Proof]{\textbf{#1.} }{\hfill$\Box$}

\def\flow{\phi_{\mathrm{low}}}

\def\R{\mathbb{R}}
\def\N{\mathbb{N}}
\def\1st{\mathrm{g}}
\def\2nd{\mathrm{H}}

\def\epsf{\epsilon_{f}}
\def\epsg{\epsilon_{\1st}}

\def\T{\mathrm{T}}
\def\Ir{\mathcal{R}}
\def\It{\mathcal{N}}
\def\cIt{c_{\It}}
\def\cIr{c_{\Ir}}
\def\sigmad{\sigma_{\mathrm{d}}}
\def\sigmaphi{\sigma_{\phi}}
\def\kappad{\kappa_{\mathrm{d}}}
\def\pq{\tfrac{1+p}{2}}
\def\alphapq{\revised{\theta_{p}}}
\def\alphap1{\revised{\theta_{1}}}

\def\halphapq{\revised{\hat{\alpha}_p}}
\def\halphap1{\revised{\hat{\alpha}_{1}}}

\numberwithin{equation}{section}

\newcommand{\revised}[1]{\textcolor{black}{#1}}
\newcommand{\morevised}[1]{\textcolor{red}{#1}}

\begin{document}

\title{A line search framework with restarting for noisy optimization
problems} 
\author{
	Albert S. Berahas
	\thanks{Department of Industrial \& Operations Engineering, 
	University of Michigan, 
	Ann Arbor, MI 48109-2102, USA
	({\tt albertberahas@gmail.com}).}
\and
	Michael J. O'Neill
	\thanks{Department of Statistics and Operations Research, 
	University of North Carolina at Chapel Hill, 
	Chapel Hill, NC 27514, USA 
	({\tt mikeoneill@unc.edu}).}
\and 
	Cl\'ement W. Royer
	\thanks{LAMSADE, CNRS, 
	Universit\'e Paris Dauphine-PSL, 75016 Paris, France
	({\tt clement.royer@lamsade.dauphine.fr}).}
}
\maketitle
\thispagestyle{empty}
{\small
\begin{abstract}
Nonlinear optimization methods are typically iterative and make use of gradient 
information to determine a direction of improvement and function information 
to effectively check for progress. When this information is corrupted by noise, 
designing a convergent and practical algorithmic process becomes challenging, 
as care must be taken to avoid taking bad steps due to erroneous information. 
For this reason, simple gradient-based schemes have been quite popular, despite 
being outperformed by more advanced techniques in the noiseless setting. 
In this paper, we propose a general algorithmic framework based on line  
search that is endowed with iteration and evaluation complexity guarantees 
even in a noisy setting. These guarantees are obtained as a result
of a restarting condition, that monitors desirable properties for the steps 
taken at each iteration and can be checked even in the presence of noise. 
Experiments using a nonlinear conjugate gradient variant and a quasi-Newton 
variant illustrate that restarting can be performed without compromising 
practical efficiency and robustness.
\end{abstract} 
}

\section{Introduction}
\label{sec:intro}

Gradient descent techniques are iterative, simple, require only first-order 
(gradient) information, and have emerged as some of the most versatile 
algorithms for smooth nonlinear optimization. In particular, their application 
in the presence of noise has received significant attention from the 
optimization and machine learning 
communities~\cite{LBottou_FECurtis_JNocedal_2018}. Under 
appropriate assumptions on the noise, complexity guarantees can be derived for 
gradient-type methods, even when those techniques use a line search in 
lieu of fixed or predetermined step size sequences~\cite{ASBerahas_LCao_KScheinberg_2021,BJin_KScheinberg_MXie_2024}. 
That said, there exist numerous advanced optimization routines that are known 
to outperform gradient descent in practice for noiseless problems, such as 
quasi-Newton and nonlinear conjugate gradient 
methods~\cite{JNocedal_SJWright_2006}. These algorithms are notoriously 
difficult to endow with complexity guarantees in the nonconvex (or sometimes 
even convex) setting. In fact, it is unclear whether those methods can exhibit 
better theoretical guarantees than simple gradient descent schemes, in a way 
that matches their practical performance. The situation is even more 
complicated in settings where the function and gradient information these 
procedures rely upon to generate and validate search directions is corrupted 
with noise.

In spite of these challenges, complexity results have recently been obtained 
by considering problems for which access to function and gradient 
information is corrupted with noise of a known or estimable level. Problems of 
this type naturally arise in numerous fields such as machine 
learning~\cite{JCDuchi_MIJordan_MJWainwright_AWibisono_2015,
MFazel_RGe_SKakade_MMesbahi_2018}, 
black-box optimization~\cite{CAudet_WHare_2017,ARConn_KScheinberg_LNVicents_2009, 
JLarson_MMenickelley_SMWild_2019}, 
and simulation optimization~\cite{RPasupathy_PGlynn_SGhosh_FSHashemi_2018,
SShashaani_FSHashemi_RPasupathy_2018}. Several algorithms based on globalization 
techniques (e.g., line/step search, trust-region) \revised{and quasi-Newton 
updates have been proposed in recent years to tackle this class of 
problems~\cite{ASBerahas_RHByrd_JNocedal_2019,
EBerglund_JZhang_MJohansson_2025,BIrwin_EHaber_2023,
xie2020analysis,shi2022noise}.} In a parallel line of work, various algorithmic 
strategies have been developed for the regime where the gradient estimates are 
sufficiently accurate with reasonable 
probability~\cite{ASBerahas_LCao_KScheinberg_2021,CCartis_KScheinberg_2018,
BJin_KScheinberg_MXie_2024,CPaquette_KScheinberg_2020,cao2024first,sun2023trust}. 
While complexity results are available for all of the aforementioned methods, 
almost all of them are based on simple, gradient descent-type steps.

In this work, we propose a general algorithmic framework of line search type for 
the unconstrained minimization of noisy functions endowed with worst-case 
complexity guarantees. We build upon the restarting conditions developed for 
deterministic nonlinear conjugate gradient methods by Chan--Renous-Legoubin and 
Royer~\cite{RChanRenousLegoubin_CWRoyer_2022} and prove that these conditions 
remain sufficient for establishing worst-case complexity results even in the 
presence of noise, under a noise condition that is on par with existing 
literature. Our algorithmic framework adaptively chooses the step size 
at each iteration via a robust line search procedure that uses the estimated 
noise level present in function evaluations. We do not place restrictions on 
how search directions are generated in our framework, enabling it to be used in 
conjunction with a wide range of popular nonlinear optimization procedures 
including quasi-Newton and nonlinear conjugate gradient.

The paper is organized as follows. We describe our problem setting and 
algorithmic framework based on restarts in Section~\ref{sec:pbalgo}. We then 
derive complexity results for our framework in Section~\ref{sec:wcc}, where we 
also formalize our assumptions regarding the noise. In Section~\ref{sec:num}, 
we demonstrate the effectiveness of our approach numerically by comparing 
implementations of our framework based on quasi-Newton and nonlinear conjugate 
gradient methods. Finally, we provide some concluding remarks in 
Section~\ref{sec:conc}.

\section{Problem and algorithm}
\label{sec:pbalgo}

This paper focuses on the problem
\begin{equation}
\label{eq:pb}
	\min_{x \in \R^n} \phi(x),
\end{equation}
where $\phi: \R^n \rightarrow \R$ is a nonconvex, continuously differentiable 
function. In our setting, the objective function and its associated derivatives 
cannot be computed exactly. Instead, approximations that are contaminated with 
noise are available. We denote the function and gradient approximations of 
$\phi$ as $f$ and $g$, respectively. The precise conditions on the quality of 
these approximations are stated in Section~\ref{sec:wcc}. Given a tolerance 
$\epsilon \in (0,1)$, our goal is to compute an $\epsilon$-approximate 
stationary point, that is, a vector $x \in \R^n$ such that
\begin{equation}
\label{eq:epspoint}
	\|\nabla \phi(x)\| \le \epsilon,
\end{equation}
where $\|\cdot\|$ denotes the Euclidean norm in $\R^n$. 

In this section, we extend the nonlinear CG method proposed by 
Chan--Renous-Legoubin and Royer~\cite{RChanRenousLegoubin_CWRoyer_2022} to
allow for inexact function and gradient values and 
describe a nonlinear conjugate gradient method based on 
Armijo line search and a modified restart condition.

\subsection{Algorithmic framework}
\label{subsec:algo}

Our proposed algorithmic framework is given in Algorithm~\ref{alg:algo}. At 
every iteration, we perform a backtracking line search to compute a step that 
yields a suitable decrease in the objective function (see 
condition~\eqref{eq:lscond}). 
\revised{Note that this sufficient decrease condition \eqref{eq:lscond} is 
relaxed by a parameter $\epsf$ in order to account for potential noise in the 
objective function evaluations. This parameter is explicitly defined in 
Assumption~\ref{as:accfun}.}

\begin{algorithm}[ht!]
\caption{Line search framework with restart conditions and noisy estimates}
\label{alg:algo}
\begin{algorithmic}[1]
    \item[] \textbf{Inputs:} $x_0 \in \mathbb{R}^d$, 
    $\eta \in \left(0,\tfrac{1}{2}\right]$, 
	$\revised{\rho} \in (0,1)$, $\sigmad \in (0,1]$, 
	$\kappad \ge 1$, $p \ge 0$\revised{, $\epsf \geq 0$}.
	\STATE Compute an estimate $g_0$ of $\nabla \phi(x_0)$.
	\STATE Set $d_0=-g_0$ and $k = 0$.
	\FOR{$k=0,1,2,\dots$}
  		\STATE Compute $\alpha_k= \revised{\rho}^{j_k}$ where $j_k$ is the smallest 
  		nonnegative integer such that
  		\begin{equation}
  		\label{eq:lscond}
  			f(x_k+\alpha_k d_k) < f(x_k) + \eta \alpha_k g_k^T d_k 
  			+ 2 \epsf.
  		\end{equation}
  		where $f(x_k+\alpha_k d_k)$ and $f(x_k)$ are estimates of 
  		$\phi(x_k+\alpha_k d_k)$ and $\phi(x_k)$, respectively.
  		\STATE Set $x_{k+1} = x_k + \alpha_k d_k$ and compute an estimate 
        $g_{k+1}$ of the gradient $\nabla \phi(x_{k+1})$.
        \STATE Define a direction $d_{k+1}$ using $g_{k+1}$.
        \STATE If the condition
  		\begin{equation}
  		\label{eq:restartcond}
  			g_{k+1}^\T d_{k+1} \ge -\sigmad\,\|g_{k+1}\|^{1+p} 
  			\quad \mbox{or} \quad 
  			\|d_{k+1}\| \ge \kappad\,\|g_{k+1}\|^{\pq},
  		\end{equation}
  		holds, restart the algorithm by setting $d_{k+1}=-g_{k+1}$.
  	\ENDFOR
\end{algorithmic}
\end{algorithm}

Once the new point has been computed, we estimate the gradient at the new 
iterate, then use that estimate to compute a new search direction.
The key ingredient in Algorithm~\ref{alg:algo} is the restarting 
condition~\eqref{eq:restartcond}, that determines whether the search direction 
(e.g., a nonlinear CG direction or a quasi-Newton direction)
is kept for the next iteration.
For any $k \ge 1$, if iteration $k-1$ does not end with a restart, it follows 
that\footnote{Although the inequalities should be strict, we use non-strict 
inequalities for notational convenience.}
\begin{equation}
\label{eq:conddk}
	g_k^\T d_k \le -\sigmad \|g_k\|^{1+p}
	\quad \mbox{and} \quad 
	\|d_k\| \le \kappad \|g_k\|^{\pq},
\end{equation}
a condition more general than that commonly used in gradient-type 
methods for both deterministic and stochastic optimization~\cite{
ASBerahas_LCao_KScheinberg_2021,
CCartis_PhRSampaio_PhLToint_2015,CCartis_KScheinberg_2018,LGrippo_SLucidi_2005}. 
In this paper, we \revised{use $\pq$ in the second part of~\eqref{eq:conddk} rather 
than introducing an additional parameter.} In addition to simplifying the 
notation, this condition was shown to be beneficial in terms of line search 
iterations in the deterministic setting~\cite{RChanRenousLegoubin_CWRoyer_2022}.

On the other hand, if iteration $k-1$ ends with a 
restart~\eqref{eq:restartcond}, it follows that $d_k=-g_k$, in which 
case~\eqref{eq:conddk} is satisfied with \revised{$p=\kappad=\sigmad=1$}. 
More generally, when \revised{$p=1$}, condition~\eqref{eq:conddk} becomes
\begin{equation}
\label{eq:conddk11}
	g_k^\T d_k \le -\sigmad\|g_k\|^2 
	\quad \mbox{and} \quad 
	\|d_k\| \le \kappad \|g_k\|.
\end{equation}
Such a condition is typical of gradient-related directions, 
and has been instrumental in obtaining complexity guarantees for 
gradient-type methods in the noiseless 
setting~\cite{CCartis_PhRSampaio_PhLToint_2015,
RChanRenousLegoubin_CWRoyer_2022}. In the 
context of nonlinear conjugate gradient methods, similar properties 
have been used to establish global convergence~\cite{WWHager_HZhang_2005}.

\subsection{Algorithmic instances}
\label{ssec:exalgos}

The framework described in Algorithm~\ref{alg:algo} is quite generic, and 
covers a number of search direction choices. In our experiments, we will 
focus on two variants corresponding to classical nonlinear optimization 
methods.

\paragraph{Nonlinear conjugate gradient} Nonlinear conjugate gradient (CG) 
methods set $d_0=-g_0$ and $d_{k+1} = -g_{k+1}+\beta_{k+1}\revised{d_k}$, 
where $\beta_{k+1}$ is a conjugate parameter~\cite{WWHager_HZhang_2006a}. 
These methods can be classified as momentum-based methods, in that they 
combine current gradient information with the previously obtained 
search direction~\cite{SJWright_BRecht_2022}.

\paragraph{Quasi-Newton methods} Quasi-Newton algorithms set $d_k = -H_k g_k$, 
where $H_k$ is a symmetric positive definite approximation of the Hessian 
matrix. Quasi-Newton methods such as (L-)BFGS typically include a restart 
condition to avoid updates that would lead to indefiniteness of 
$H_{k+1}$~\cite[Chapter 6]{JNocedal_SJWright_2006}. More precisely, the update 
is skipped (i.e. $H_{k+1} = H_k$) if
\begin{equation}
\label{eq:restartbfgs}
	s_k^\T (g_{k+1}-g_k) \ge \tau \|s_k\| \|g_{k+1}-g_k\|,
\end{equation}
fails, where $s_k=\alpha_k d_k=x_{k+1}-x_k$ and $\tau \in (0,1)$. This 
process guarantees that $H_{k+1} \succ 0$. Interestingly, 
when~\eqref{eq:restartbfgs} holds, we also have
\[
	g_{k+1}^\T (-H_{k+1} g_{k+1}) \le -\lambda_{\min}(H_{k+1})\|g_{k+1}\|^2
	\quad \mbox{and} \quad
	\|H_{k+1}g_{k+1}\| \le \lambda_{\max}(H_{k+1})\|g_{k+1}\|,
\]
\revised{and those conditions resemble~\eqref{eq:conddk11} with $p=1$.}

\section{Complexity analysis}
\label{sec:wcc}

In this section, we derive a complexity result for our line search framework 
with restart conditions and noisy estimates. Section~\ref{subsec:wcc:asslem} 
provides the necessary assumptions as well as intermediate results, while 
Section~\ref{subsec:wcc:main} establishes and discusses the complexity bound 
for our algorithm.

\subsection{Assumptions and decrease guarantees}
\label{subsec:wcc:asslem}

We make the following assumptions about the objective function of 
problem~\eqref{eq:pb}.

\begin{assumption}
\label{as:fC11L}
	The function $\phi$ is continuously differentiable on $\R^n$ and its 
	gradient is $L$-Lipschitz continuous for $L>0$.
\end{assumption}

\begin{assumption}
\label{as:flow}
	There exists $\flow \in \R$ such that $\phi(x) \ge \flow$ for every 
	$x \in \R^n$.
\end{assumption}

\revised{We now state an assumption regarding the objective function 
estimates, that involves the parameter $\epsf$ used in 
Algorithm~\ref{alg:algo}.}

\begin{assumption}
\label{as:accfun}
    There exists a constant $\epsf \ge 0$ such that for all $x \in \R^n$ 
    the objective function estimate $f(x)$ of $\phi(x)$
	\begin{equation}
	\label{eq:accfun}
		| f(x) - \phi(x) | \le \epsf.
	\end{equation}
\end{assumption}

The next assumption pertains to the objective function gradient estimates. 
This assumption is only required to hold for the iterates generated by 
Algorithm~\ref{alg:algo}.
\begin{assumption}
\label{as:accgrad}
    There exist constants $\sigma_{\phi} \in [0,1)$ and 
    $\epsg \ge 0$ such that for all $k \in \N$ the objective gradient 
    estimate at iteration $k$ (denoted by $g_k$) satisfies
    \begin{equation}
	\label{eq:accgrad}
        \|g_k - \nabla \phi(x_k)\|
        \le 
		\max\left\{\epsg,\sigmaphi\alphapq 
		\min\left[\|\nabla \phi(x_k)\|,
		\|\nabla \phi(x_k)\|^{\pq}  \right] \right\}
    \end{equation}
    where $x_k \in \mathbb{R}^n$ are the iterates generated by 
    Algorithm~\ref{alg:algo}, $p \geq 0$ is the value used 
    in~\eqref{eq:restartcond}, and $\alphapq$ is defined as
    \begin{equation}
    \label{eq:baralpha}
    	\alphapq \revised{:= 
        \frac{(1-\eta)\sigmad(1-\sigmaphi)^{1+p}}{2\kappad(1+\sigmaphi)^{\pq}}.}
    \end{equation}
\end{assumption}

Assumption~\ref{as:accgrad} resembles conditions previously used in the 
literature in that it generalizes norm conditions~\cite{carter1991global} 
as well as bounded noise~\cite{ASBerahas_RHByrd_JNocedal_2019}. Although we 
consider the true gradient on the right-hand side of 
condition~\eqref{eq:accgrad}, we point out that the gradient estimate could 
also be used, akin to~\cite{ASBerahas_LCao_KScheinberg_2021,byrd2012sample,
CCartis_KScheinberg_2018}.
Note that when \revised{$p=1$}, condition~\eqref{eq:accgrad} becomes
\begin{align*}
    \|g_k - \nabla \phi(x_k)\|
    \le    
    \max\left\{\epsg,\sigmaphi\alphapq\|\nabla \phi(x_k)\|
    \right\}.
\end{align*}
When $\sigmaphi=0$, condition~\eqref{eq:accgrad} corresponds to bounded 
noise in the gradient estimates. Conversely, when $\epsg=0$, the condition 
becomes a norm condition involving the true gradient norm as well as a 
constant coefficient~\cite{carter1991global}. We note that some conditions 
proposed in the literature involve $\alpha_k \|g_k\|$ on the right-hand 
side rather than a constant  times 
$\|\nabla \phi(x_k)\|$~\cite{ASBerahas_LCao_KScheinberg_2021}. 
\revised{We focus on using the true gradient norm $\|\nabla \phi(x_k)\|$ 
for simplicity in this paper, but we use 
a fixed quantity $\alphapq$ in lieu of the classical iteration-dependent 
quantity $\alpha_k$. The former ($\alphapq$)} \revised{is sufficient to derive} {a 
lower bound on the step size,  which typically exists for line search 
schemes.} \revised{Using $\alphapq$, which is associated with our worst-case lower-bound on $\alpha_k$, instead of $\alpha_k$ may represent a stronger condition on some iterations. However, for sufficiently small $\alpha_k$, this may not be the case, as $\alphapq$ is \textit{independent} of the Lipschitz constant $L$, while any worst-case lower bound on $\alpha_k$ will be dependent on this quantity.}

More broadly, our focus in this paper is to examine the effect 
of the restarting conditions in the noisy setting, and thus we focused 
on a single noise condition for that purpose.
%
%
As we will see below, the use of $\|\nabla \phi(x_k)\|^{\pq}$ combined 
with~\eqref{eq:conddk} is instrumental for deriving complexity results.
\revised{Those results are a generalization of the noiseless 
setting~\cite{RChanRenousLegoubin_CWRoyer_2022}, that can be recovered 
for $\epsilon_f=\epsilon_g=0$.}

\begin{lemma}
\label{le:ineqngk}
    Let Assumption~\ref{as:accgrad} hold, and suppose that at the $k$th 
    iteration of Algorithm~\ref{alg:algo}
    \begin{equation}
    \label{eq:gradbigenough}
        \min\left\{\|\nabla \phi(x_k)\|,\|\nabla \phi(x_k)\|^{\pq} \right\} 
    	\ge 
    	\frac{\epsg}{\sigmaphi\alphapq},
    \end{equation}
    where $\alphapq$ is defined \revised{in~\eqref{eq:baralpha}}. Then, 
    \begin{equation}
    \label{eq:gklowerbound}
		\|g_k\| 
		\ge 
		\left(1-\sigmaphi\right)\|\nabla \phi(x_k)\|
    \end{equation}
    and 
    \begin{equation}
    \label{eq:gkupperbound}
        \|g_k\| \le (1+\sigmaphi)\|\nabla \phi(x_k)\|.
    \end{equation}
\end{lemma}

\begin{proof}
We begin by proving~\eqref{eq:gklowerbound}. By Assumption~\ref{as:accgrad}, 
we have
\begin{eqnarray}
\label{eq:auxlbgk}
	\|\nabla \phi(x_k)\|
	&\le 
	&\|\nabla \phi(x_k)-g_k\| + \|g_k\| \nonumber \\
	&\le 
	&\max\left\{
	\epsg,\sigmaphi\,\alphapq\,\min\left\{\|\nabla \phi(x_k)\|,
	\|\nabla \phi(x_k)\|^{\pq} \right\} \right\} +\|g_k\| \nonumber \\
    &\le 
    &\sigmaphi\,\alphapq\,\min\left\{\|\nabla \phi(x_k)\|,
	\|\nabla \phi(x_k)\|^{\pq} \right\}  + \|g_k\| \nonumber\\
    &\le 
    &\sigmaphi\,\|\nabla \phi(x_k)\| + \|g_k\|,
\end{eqnarray}
where the third inequality follows from~\eqref{eq:gradbigenough} and 
the last inequality uses $\alphapq \le 1$ \revised{(which follows 
from the definition~\eqref{eq:baralpha} along with 
$\eta \in (0,\tfrac{1}{2}]$, $\sigmaphi \in [0,1)$, $\sigmad \in (0,1]$ and 
$\kappad \ge 1$)}.  
Re-arranging the terms in the 
last inequality yields~\eqref{eq:gklowerbound}.

We now turn to
~\eqref{eq:gkupperbound}. Using again the triangle 
inequality together with~\eqref{eq:gradbigenough} and $\alphapq \le 1$ 
yields
\begin{eqnarray}
\label{eq:auxubgk}
    \|g_k\| 
    &\le &\|\nabla \phi(x_k)-g_k\| + \|\nabla \phi(x_k)\| 
    \nonumber \\
    &\le &\max\left\{
    \epsg,\sigmaphi\alphapq\min\left\{\|\nabla \phi(x_k)\|,
    \|\nabla \phi(x_k)\|^{\pq} \right\}
    \right\} + \|\nabla \phi(x_k)\|  
    \nonumber \\
    &\le &\sigmaphi\alphapq \min\left\{\|\nabla \phi(x_k)\|,
    \|\nabla \phi(x_k)\|^{\pq} \right\} + \|\nabla \phi(x_k)\| 
    \nonumber \\
    &\le &\sigmaphi\min\left\{\|\nabla \phi(x_k)\|,
    \|\nabla \phi(x_k)\|^{\pq} \right\} + \|\nabla \phi(x_k)\| \\
    &\le &\left(1+\sigmaphi\right)\|\nabla \phi(x_k)\|,
\end{eqnarray}
and thus~\eqref{eq:gkupperbound} holds.
\end{proof}

To establish a complexity result for our algorithmic framework 
we divide the 
iterations into two categories, identified by the indices
\begin{equation}
\label{eq:wcc:itpart}
	\begin{array}{lll}
		\Ir &= 
		&\{ k \in \N\ |\ 
		d_k = -g_k \} \\
		\It &= &\N \setminus \Ir.
	\end{array}
\end{equation}
Any $k \in \Ir$ is the index of a \emph{restarted iteration}, while 
$k \in \It$ is the index of a \emph{non-restarted iteration}.
Depending on the nature of each iteration, we can bound the number 
of backtracking steps needed to compute a suitable step size. We begin with 
the non-restarted iterations, as the proof encompasses that of restarted 
iterations.

\begin{lemma}
\label{le:wcc:lst}
	Let Assumptions~\ref{as:fC11L}, \ref{as:accfun}, and \ref{as:accgrad} 
	hold. Suppose that the $k$th iteration of Algorithm~\ref{alg:algo} is 
    a non-restarted iteration (i.e., $k \in \It$), and that
    \begin{equation}
    \label{eq:nonstatiolst}
        \min\left\{\|\nabla \phi(x_k)\|,\|\nabla \phi(x_k)\|^{\pq}\right\} 
        \ge \frac{\epsg}{\sigmaphi\alphapq},
    \end{equation}
    where $\alphapq$ is defined in~\eqref{eq:accgrad}.
	Then, the line search process terminates after at most 
	$\lfloor \bar{j}_{\It}+1 \rfloor$ iterations, where 
	\begin{equation}
	\label{eq:itslt}
		\bar{j}_{\It} := 
		\left[ \log_{\revised{\rho}}(\alphapq) 
		\right]_+     
	\end{equation}
	\revised{and $[t]_+ = \max\{t,0\}$.}
	Moreover, the resulting decrease at the $k$th iteration satisfies
	\begin{equation}
	\label{eq:declt}
		\phi(x_k) - \phi(x_{k+1}) 
		>
        \cIt \epsg^{1+p}-4\epsf
	\end{equation}
	where 
	\[
		c_{\It}:= \eta\sigmad\frac{(1-\sigmaphi)^{1+p}}{\sigmaphi^{\revised{1+p}}\revised{\theta_p^{1+p}}}
		\min\{1,\revised{\rho\halphapq}\}, \qquad \revised{and} \qquad \halphapq := \frac{(1-\eta)\sigmad(1-\sigmaphi)^{1+p}}{L\kappad^2(1+\sigmaphi)^{1+p}}.
	\]
\end{lemma}

\begin{proof}
	Suppose first that $\alpha_k = 1$ satisfies the Armijo-type condition. 
	In that case, ~\eqref{eq:itslt} holds, and 
	condition~\eqref{eq:lscond} together with $k \in \It$ gives
	\begin{eqnarray}
	\label{eq:declsIt}
		f(x_k)-f(x_{k+1}) > -\eta g_k^\T d_k - 2 \epsf 
		> \eta \sigmad \|g_k\|^{1+p} - 2 \epsf.
	\end{eqnarray}
    Combining Assumption~\ref{as:accfun}, \eqref{eq:gklowerbound}, 
    \eqref{eq:nonstatiolst} and \eqref{eq:declsIt} together with the fact 
    that $\sigmaphi\alphapq\le 1$, it follows that	
	\begin{eqnarray*}
		\phi(x_k)-\phi(x_{k+1}) 
		&\ge 
		&f(x_k)-f(x_{k+1}) - 2\epsf \\
		&\ge 
		&\eta\sigmad \|g_k\|^{1+p} - 4 \epsf \\
		&\ge 
		&\eta\sigmad (1-\sigmaphi)^{1+p}\|\nabla \phi(x_k)\|^{1+p} 
		- 4 \epsf \\
		&\ge 
		&\eta\sigmad (1-\sigmaphi)^{1+p} \left[
		\frac{\epsg}{\sigmaphi\alphapq}\right]^{1+p}
		- 4\epsf \\
		&\ge 
		&\cIt \epsg^{1+p}-4\epsf.
	\end{eqnarray*}
	As a result, \eqref{eq:declt} holds in that case.

    Suppose now that the line search condition fails for some 
	$\alpha=\revised{\rho}^j$ with $j \in \N$, i.e., \eqref{eq:lscond} does not hold,
	\begin{equation*}
		\eta \alpha g_k^\T d_k 
		\le f(x_k+\alpha d_k)-f(x_k) - 2 \epsf
	\end{equation*}
    By Assumption~\ref{as:fC11L} and \eqref{eq:accfun}, it follows that
    \begin{align*}
		\eta \alpha g_k^\T d_k 
		&\le f(x_k+\alpha d_k)-f(x_k) - 2 \epsf\\
            &\le \phi(x_k+\alpha d_k)-\phi(x_k)\\
            &\le \alpha \nabla \phi(x_k)^\T d_k + \frac{L}{2}\alpha^2 \|d_k\|^2.
	\end{align*}
    Subtracting $g_k^\T d_k$ from both sides and dividing by $\alpha$
    \begin{align*}
		-(1-\eta)g_k^\T d_k 
		&\le \left[\nabla \phi(x_k)-g_k\right]^\T d_k + \frac{L}{2}\alpha \|d_k\|^2 \\
		&\le \|\nabla \phi(x_k)-g_k\| \|d_k\|  + \frac{L}{2}\alpha \|d_k\|^2.
	\end{align*}
    Since $k \in \It$, condition~\eqref{eq:conddk} is satisfied. Using 
    both inequalities in \eqref{eq:conddk}, it follows that the above 
    inequality can be bounded as follows
	\begin{eqnarray*}
		(1-\eta)\sigmad \|g_k\|^{1+p} 
		&\le &\kappad \|\nabla \phi(x_k)-g_k\| \|g_k\|^{\pq} 
		+ \frac{L\kappad^2}{2}\alpha \|g_k\|^{(1+p)}.
	\end{eqnarray*}
    \revised{By~\eqref{eq:accgrad} and~\eqref{eq:nonstatiolst},}
	\begin{align*}
        (1-\eta)\sigmad \|g_k\|^{1+p}
        &\le 
        \kappad \sigmaphi\alphapq \min\left\{\|\nabla \phi(x_k)\|,
        \|\nabla \phi(x_k)\|^{\pq} \right\} \|g_k\|^{\pq} 
        + \frac{L\kappad^2}{2}\alpha \|g_k\|^{1+p} \nonumber \\
        &\le 
        \kappad \sigmaphi \alphapq \|\nabla \phi(x_k)\|^{\pq} \|g_k\|^{\pq}
        + \frac{L\kappad^2}{2}\alpha \|g_k\|^{1+p}. \nonumber  
	\end{align*}
    Using Lemma~\ref{le:ineqngk} to bound the $\|g_k\|$ terms yields
	\begin{eqnarray}
    \label{eq:gkfailedcond}
		& 
		&(1-\eta)\sigmad(1-\sigmaphi)^{1+p} \|\nabla \phi(x_k)\|^{1+p}
		 \nonumber \\
        &\le 
        &\kappad\sigmaphi(1+\sigmaphi)^{\pq} \alphapq\|\nabla\phi(x_k)\|^{1+p} 
        + \frac{L\kappad^2(1+\sigmaphi)^{1+p}}{2}\alpha \|\nabla \phi(x_k)\|^{1+p}.
	\end{eqnarray}

    \revised{By the definition of $\alphapq$ and $\sigmaphi \leq 1$, we have
    \begin{eqnarray*}
        \kappad\sigmaphi(1+\sigmaphi)^{\pq} \alphapq\|\nabla\phi(x_k)\|^{1+p} 
        \leq 
        \frac12(1-\eta)\sigmad(1-\sigmaphi)^{1+p} \|\nabla \phi(x_k)\|^{1+p},
    \end{eqnarray*}
    and thus,
    \begin{equation*}
		\frac{(1-\eta)\sigmad(1-\sigmaphi)^{1+p}}{2} \|\nabla \phi(x_k)\|^{1+p}
        \le \frac{L\kappad^2(1+\sigmaphi)^{1+p}}{2}\alpha \|\nabla \phi(x_k)\|^{1+p}.
	\end{equation*}
    Since $\|\nabla \phi(x_k)\|>0$ by~\eqref{eq:nonstatiolst}, this inequality can 
	be further rewritten as
    \begin{equation*}
        \alpha \geq \frac{(1-\eta)\sigmad(1-\sigmaphi)^{1+p}}{L\kappad^2(1+\sigmaphi)^{1+p}} = \halphapq.
    \end{equation*}}
    
    Overall, we have shown that if the decrease condition~\eqref{eq:lscond} fails for a step size $\alpha$, then
    \begin{equation}
    \label{eq:boundalpha}
        \alpha \ge \revised{\halphapq.}
    \end{equation}
    The condition~\eqref{eq:boundalpha} can only hold if 
    $j \le \bar{j}_{\It}$, and thus the backtracking process must 
    terminate after $j_k \le \lceil \bar{j}_{\It} \rceil$ 
    iterations. In addition, since this process did not terminate with 
    $\alpha=\revised{\rho}^{j_k-1}$, it must be that
    \[
    	\revised{\rho}^{j_k-1} \ge \revised{\halphapq},
    \]
    hence, $\alpha_k=\revised{\rho}^{j_k} \ge \revised{\rho \halphapq}$.
    By the line search condition~\eqref{eq:lscond}, the fact that 
    $k \in \It$ and \eqref{eq:nonstatiolst}, 
 	\begin{eqnarray*}
		\phi(x_k)-\phi(x_{k+1}) 
		&\ge &f(x_k)-f(x_{k+1}) - 2\epsf \\
		&\ge &-\eta\alpha_k g_k^\T d_k - 4 \epsf \\
		&\ge &\eta\revised{\rho \halphapq}\sigmad 
		\|g_k\|^{1+p} - 4 \epsf \\
		&\ge 
		&\eta\revised{\rho \halphapq}\sigmad(1-\sigmaphi)^{1+p}
		\|\nabla \phi(x_k)\|^{1+p} - 4 \epsf \\
		&\ge 
		&\eta\revised{\rho \halphapq}\sigmad (1-\sigmaphi)^{1+p} \left[
		\frac{\epsg}{\sigmaphi\alphapq}\right]^{1+p}
		-4\epsf \\
		&\ge &\cIt \epsg^{1+p}-4\epsf,
	\end{eqnarray*}   
	and hence~\eqref{eq:declt} also holds in the latter case. 
    Combining both cases finally yields the desired result.
\end{proof}

We now consider the restarted iterations, i.e., $k \in \Ir$. 
In this case, we have $d_k=-g_k$.

\begin{lemma}
\label{le:wcc:lsr}
	Let Assumptions~\ref{as:fC11L}, \ref{as:accfun} and \ref{as:accgrad} 
	hold. Suppose that the $k$th iteration of Algorithm~\ref{alg:algo} is a 
	restarted iteration (i.e., $k \in \Ir$), and that
    \begin{equation}
    \label{eq:nonstatiolsr}
        \|\nabla \phi(x_k)\| \ge \frac{\epsg}{\sigmaphi\alphap1},
        \qquad \revised{where} \qquad 
        \revised{\alphap1=
        \frac{(1-\eta)(1-\sigmaphi)^2}{2(1+\sigmaphi)}.}
    \end{equation}
    Then, the line search process terminates after at most 
	$\lfloor \bar{j}_{\Ir}+1 \rfloor$ iterations, where 
	\begin{equation}
	\label{eq:itslr}
		\bar{j}_{\Ir} := 
		\left[ \log_{\revised{\rho}}(\alphap1) 
		\right]_+.
	\end{equation}
	Moreover, the resulting decrease at the $k$th iteration satisfies
	\begin{equation}
	\label{eq:declr}
		\phi(x_k) - \phi(x_{k+1}) 
		>
		\cIr \epsg^2-4\epsf
	\end{equation}
	where 
	\[
		c_{\Ir}:= \revised{\eta\frac{(1-\sigmaphi)^2}{\sigmaphi^2\alphap1^2}}\min\{1,\revised{\rho\halphap1}\} \qquad \revised{and \qquad
        \halphap1 := \frac{(1-\eta)(1-\sigmaphi)^2}{L(1+\sigmaphi)^2}.}
	\]
\end{lemma}

\begin{proof}
	The proof is a mere restatement of that of Lemma~\ref{le:wcc:lst} 
	\revised{with special values for the parameters. Indeed, at every restarted 
	iteration, condition~\eqref{eq:conddk} holds with}
	\begin{align*}
		\revised{\sigmad=1,\ \kappad=1,\ p=1},
	\end{align*}
	so that $\alphapq=\alphap1$. \revised{The rest of the proof follows.}
\end{proof}

The results of Lemmas~\ref{le:wcc:lst} and~\ref{le:wcc:lsr} are instrumental 
to bounding the number of iterations necessary to reach an approximate
stationary point.

\subsection{Main results}
\label{subsec:wcc:main}

Our main result is a bound on the number of iterations required by the 
algorithm prior to reaching an approximate stationary point~\eqref{eq:epspoint}. 
This bound also applies to the number of objective function gradient evaluations.

\begin{theorem}
\label{th:wcc:gdits}
	Let Assumptions~\ref{as:fC11L} and~\ref{as:flow} hold. 
	Suppose further that Assumptions~\ref{as:accfun} and~\ref{as:accgrad} hold 
	with 
	\begin{equation}
	\label{eq:wcc:epsgf}
		\epsilon_f 
		\; \le \;
		\min\left[ 
        \frac{c_{\It}}{8}\epsg^{1+p},
		\frac{c_{\Ir}}{8}\epsg^2
		\right]
	\end{equation}
	and let 
	\[
		\epsilon := \max\left\{\tfrac{\epsg}{\sigmaphi\alphapq},
        \left[\tfrac{\epsg}{\sigmaphi\alphapq}\right]^{\tfrac{2}{1+p}},
		\tfrac{\epsg}{\sigmaphi\alphap1}\right\}.
	\]
	Then, the number of 
	iterations (and gradient estimates) required by 
	Algorithm~\ref{alg:algo} to reach a point satisfying~\eqref{eq:epspoint} 
	is at most
	\begin{equation}
	\label{eq:wcc:gdits}
		K_{\epsilon}:= 
		\left\lceil \frac{2(\phi(x_0)-\flow)}{c_{\Ir}} 
		\epsg^{-2} 
        + \frac{2(\phi(x_0)-\flow)}{c_{\It}}\epsg^{-(1+p)} 
		\right\rceil.
	\end{equation}
\end{theorem}

\begin{proof}
	Let $K \in N$ be such that $\|\nabla \phi(x_k)\| > \epsilon$ for 
	any $k=0,\dots,K-1$. Following our partitioning~\eqref{eq:wcc:itpart}, 
	we define the index sets 
	\begin{eqnarray*}
		\It_K &:= &\It \cap \{0,\dots,K-1\},\\
		\Ir_K &:= &\Ir \cap \{0,\dots,K-1\}.
	\end{eqnarray*}
	
	For any $k \in \It_K$, the result of
	Lemma~\ref{le:wcc:lst} applies, and we have
	\begin{equation}
	\label{eq:dect}
		\phi(x_k) - \phi(x_{k+1}) 
    	> c_{\It} \epsg^{1+p} - 4 \epsilon_f 
		\ge \frac{c_{\It}}{2} \epsg
        ^{1+p}. 
	\end{equation}
	On the other hand, if $k \in \Ir_K$, applying Lemma~\ref{le:wcc:lsr} 
	gives
	\begin{equation}
	\label{eq:decr}
		\phi(x_k)-\phi(x_{k+1}) 
		> c_{\Ir} \epsg^2 - 4 \epsilon_f 
		\ge \frac{c_{\Ir}}{2}\epsg^2.
	\end{equation}
	
	We now consider the sum of function changes over all $k \in \{0,\dots,K-1\}$. 
	By Assumption~\ref{as:flow}, we obtain
	\begin{eqnarray*}
		\phi(x_0) - \flow &\ge &\phi(x_0) - \phi(x_K) \\
		&\ge &\sum_{k=0}^{K-1} \left[ \phi(x_k) - \phi(x_{k+1})\right] \\
		&\ge &\sum_{k \in \It_K} \left[ \phi(x_k) - \phi(x_{k+1})\right]
		+ \sum_{k \in \Ir_K}  \left[ \phi(x_k) - \phi(x_{k+1})\right] \\
		&> &\sum_{k \in \It_K} \frac{c_{\It}}{2} \epsg^{1+p} 
		+ \sum_{k \in \Ir_K} \frac{c_{\Ir}}{2} \epsg^2.
	\end{eqnarray*}
	Since the right-hand side consists of two sums of positive terms, the 
	above inequality implies that
	\[
		\phi(x_0) - \flow > \sum_{k \in \It_K} \frac{c_{\It}}{2}
		\epsg^{1+p}
		\; \Leftrightarrow \; 
		|\It_K| < \frac{2(\phi(x_0)-\flow)}{c_{\It}} 
        \epsg^{-(1+p)}
	\]
	and
	\[
		\phi(x_0) - \flow > \sum_{k \in \Ir_K} \frac{c_{\Ir}}{2} \epsg^2 
		\; \Leftrightarrow \;
		|\Ir_K| < \frac{2(\phi(x_0)-\flow)}{c_{\Ir}} \epsg^{-2}.
	\]
	Using $|\It_K|+|\Ir_K|=K$ finally yields 
	\[
		K < 
		\frac{2(\phi(x_0)-\flow)}{c_{\Ir}} \epsg^{-2} 
		+ \frac{2(\phi(x_0)-\flow)}{c_{\It}}
        \epsg^{-(1+p)},
	\]	
	hence $K \le K_{\epsilon}$.
\end{proof}

By combining the result of Theorem~\ref{th:wcc:gdits} with that of 
Lemmas~\ref{le:wcc:lst} and~\ref{le:wcc:lsr}, we can also provide an 
evaluation complexity bound of Algorithm~\ref{alg:algo}.

\begin{corollary}
\label{co:wccevals}
	Under the assumptions of Theorem~\ref{th:wcc:gdits}, 
    the number of 
	function evaluations required by Algorithm~\ref{alg:algo} to reach a 
	point satisfying~\eqref{eq:epspoint} is at most 
	\begin{equation}
	\label{eq:wccevals}
		\left(\left[\log_{\revised{\rho}}\left(\min\{\alphapq,\alphap1\}\right)
		\right]_+ + 1 \right) K_{\epsilon},
	\end{equation}
	where $K_{\epsilon}$ is defined in~\eqref{eq:wcc:gdits}.
\end{corollary}

\begin{proof}
    By Lemmas~\ref{le:wcc:lst} and~\ref{le:wcc:lsr},
    any iteration requires at most 
	\[
		\max\left\{\bar{j}_{\It},\bar{j}_{\Ir}\right\} + 1
		= \left[
		\log_{\revised{\rho}}\left(\min\{\alphapq,\alphap1\}\right)
		\right]_+ + 1
	\]
	function evaluations. Combining this number with the result of 
	Theorem~\ref{th:wcc:gdits} completes the proof.
\end{proof}

\begin{remark}
    When $p=1$, our results resemble the results 
    in~\cite{ASBerahas_LCao_KScheinberg_2021}. When $\epsf=\epsg=\sigmaphi=0$, 
    our results match the results obtained by 
    Chan--Renous-Legoubin and Royer for the restarted nonlinear conjugate 
    gradient method~\cite{RChanRenousLegoubin_CWRoyer_2022}. 
\end{remark}

\section{Numerical experiments}
\label{sec:num}

In this section, we investigate the numerical behavior of several methods that 
fit into the line search algorithmic framework with restart conditions and 
estimates given in Algorithm~\ref{alg:algo}. We evaluate the performance of all 
methods on a subset of the unconstrained optimization problems from the CUTEst 
collection~\cite{NIMGould_DOrban_PhLToint_2015} in both the noise-free and 
noisy settings. The purpose of these experiments is threefold: $(1)$ to 
highlight the merits and limitations of different methods; $(2)$ to understand 
the effect of the restart conditions; and, $(3)$ to explore the influence of 
the noise. 

\subsection{Setup}
\label{ssec:setup}

We implemented our algorithmic framework and different instances of the methods 
in MATLAB R2024a. Experiments were run on the Longleaf 
cluster\footnote{\url{https://help.rc.unc.edu/longleaf-cluster/}}, a 
Linux-based computing system that contains over 1840 CPU cores across 230 
blade servers.

\paragraph{Test problems} Our test set consists of 234 unconstrained problems from 
the CUTEst collection~\cite{NIMGould_DOrban_PhLToint_2015} with dimensions 
between 1 and 1000. We used the new MATLAB implementation and interface of the 
CUTEst problems developed by Gratton and Toint~\cite{SGratton_PhLToint_2024}, and 
the default initial points provided for those problems. We scaled each problem by 
$\max\{1,\|\nabla \phi_0\|_{\infty}\}$, where $\nabla \phi_0$ is the true gradient 
at the initial point.

We considered both noise-free (deterministic) and noisy versions of all the 
problems. For each problem in the test set, we considered 5 noise levels 
(including the noise-free setting), namely 
$\epsilon_f = \{ 0, 10^{-8}, 10^{-4}, 10^{-2}, 10^{-1}\}$, and 
$\epsg = \sqrt{\epsilon_f}$. We added uniform noise to the function and gradient 
evaluations bounded by $\epsilon_f$ and $\epsg$, respectively, and assumed that 
$\epsilon_f$ was available to all methods.

\paragraph{Algorithms} We compared the following methods.
\begin{enumerate}[(1)]
    \item \textbf{Gradient descent (GD):} Algorithm~\ref{alg:algo} with 
    $d_k=-g_k$, $p=1$ and $\kappad=\sigmad^{-1}=1$.
    \item \textbf{Nonlinear conjugate gradient (NLCG):} 
    Algorithm~\ref{alg:algo} with $d_k = - g_k + \beta_k d_{k-1}$, where 
    $\beta_k$ is updated using the PRP+ formula~\cite{WWHager_HZhang_2006b} 
    and $\sigmad = 0$ and $\kappad = \infty$.
    \item \textbf{Limited-Memory BFGS (LBFGS):} The classical LBFGS 
    method~\cite[equation (7.15)]{JNocedal_SJWright_2006}. 
    Algorithm~\ref{alg:algo} with $d_k = - H_kg_k$, where $H_k$ is updated 
    using the BFGS methodology~\cite[equation (7.19)]{JNocedal_SJWright_2006} 
    and $d_k$ is computed using the LBFGS two-loop 
    recursion~\cite[Algorithm 7.4]{JNocedal_SJWright_2006}, and $\sigmad = 0$ 
    and $\kappad = \infty$. We consider a memory size $m = 10$, and, as is 
    standard practice, curvature pairs $s_k = x_{k+1}-x_k $ and 
    $y_k = g_{k+1} - g_k$ are only stored if 
    $s_k^Ty_k \geq \epsilon_{sy} \|s_k \| \|y_k\|$, where 
    $\epsilon_{sy} = 10^{-4}$~\cite{JNocedal_SJWright_2006}. 
    \item \textbf{Nonlinear conjugate gradient with restarts (NLCGr):} 
    Algorithm~\ref{alg:algo} with $d_k = - g_k + \beta_k d_{k-1}$, where 
    $\beta_k$ is updated using the PRP+ formula~\cite{WWHager_HZhang_2006b}. 
    We choose $\kappad \in (0, \infty)$ and set $\sigmad=\kappad^{-1}$.
    \item \textbf{Limited-Memory BFGS with restarts (LBFGSr):} 
    Algorithm~\ref{alg:algo} with $d_k = - H_kg_k$, where $H_k$ is updated 
    using the BFGS methodology~\cite[equation (7.19)]{JNocedal_SJWright_2006}
    and $d_k$ is computed using the LBFGS two-loop 
    recursion~\cite[Algorithm 7.4]{JNocedal_SJWright_2006}. 
    We consider a memory size $m = 10$, and a similar strategy 
    as that of LBFGS is used to store curvature pairs at every iteration.
    We choose $\kappad \in (0, \infty)$ and set $\sigmad=\kappad^{-1}$.
\end{enumerate}

For methods (2) and (3), the values of $\sigmad=0$ and $\kappad=\infty$ 
ensure that the methods only restart if the computed direction is not a 
descent direction, i.e., if $g_k ^Td_k \geq 0$. This, of course, never 
occurs for method (3). For the restarted variants (4) and (5), we consider 
multiple values for $p$ and $\kappad$ for a sensitivity analysis, using 
$\sigmad=\kappad^{-1}$ to restrict this analysis to two parameters. 
Following the previous study on restarted nonlinear conjugate gradient and 
line search techniques~\cite{YCarmon_JCDuchi_OHinder_ASidford_2017a}, we set 
the line search parameters to $\eta = 1/2$ and $\morevised{\rho} = 1/2$. We also 
consider the PRP+ update rule for $\beta_k$ since it has proven to be robust 
and competitive in practice~\cite{RChanRenousLegoubin_CWRoyer_2022}. 
Finally, all methods were run 10 times on each problem and noise instance to 
account for randomness in function/gradient estimates. 

\paragraph{Termination condition} Algorithms were terminated when a maximum 
number of $1000$ iterations was reached, or when
$\| g_k \|_{\infty} \leq \max\{2 \epsg,10^{-8}\}$. Note that this condition 
implies that wrong termination can occur whenever the noise induces a gradient 
estimate with small norm. We have encountered this situation at the initial 
point for 2.6\% of the noisy runs, and 5.6\% of the runs for the largest noise 
level. The results presented thereafter discard those runs, and consider a 
problem to be solved whenever the method computed an iterate $x_{\hat{k}}$ 
such that
\begin{align*}
    \|\nabla \phi(x_{\hat{k}})\|_{\infty} 
    \le 
    \epsg+\max\{2\epsg,10^{-8}\}.
\end{align*}
This criterion is based on~\eqref{eq:accgrad}, and reduces to 
$\|\nabla \phi(x_{\hat{k}})\|_{\infty} \le 10^{-8}$ in the absence of noise.

\paragraph{Performance and Data Profiles} We compare restarted variants of 
NLCG and LBFGS with their non-restarted counterparts as well as gradient 
descent for several noise levels. We summarize our findings using performance 
and data profiles~\cite{EDDolan_JJMore_2002,JJMore_SMWild_2009}.
We built performance and data profiles for our tests based on the number of 
gradient calls (equivalent to the number of iterations) required to achieve 
a desired tolerance.

\subsection{Results in the noiseless setting}
\label{ssec:numnoise0}

We begin our study by investigating the noiseless setting. To this end, and 
similarly to previous work on restarting 
NLCG~\cite{RChanRenousLegoubin_CWRoyer_2022}, we first compute the average 
percentage of restarts for all variants for different choices of $p$ and 
$\kappad$ (per our settings, all restarted variants use $\sigmad=\kappad^{-1}$).

\begin{table}[ht]
{\footnotesize
\centering
\begin{tabular}{l|crrrrr}
\toprule
\multicolumn{7}{c}{NLCGr} \\
\midrule
$p\backslash\kappad$ & &$10^{2}$ &$10^{3}$ &$10^{4}$ &$10^{5}$ &$\mathbf{10^{6}}$ \\
\cmidrule{1-7}
0~~~ & &77.86 &66.66 &54.42 &40.84 &26.29 \\
0.25 & &68.09 &52.38 &34.43 &14.66 &~0.95 \\
0.50 & &48.15 &21.51 &~0.92 &~0.35 &\textbf{~0.32} \\
\textbf{0.75} & &~2.01 &~0.60 &~0.40 &~0.34 &\textbf{~0.32} \\
1~~~ & &~2.20 &~0.70 &~0.42 &~0.34 &\textbf{~0.32} \\
\bottomrule
\end{tabular}
\quad 
\begin{tabular}{l|crrrrr}
\toprule
\multicolumn{7}{c}{LBFGSr} \\
 \midrule
$p\backslash\kappad$ & &$10^{2}$ &$10^{3}$ &$10^{4}$ &$10^{5}$ &$\mathbf{10^{6}}$ \\
 \cmidrule{1-7}
0~~~ & &91.38	&81.58	 &69.21	 &46.90	 &21.82 \\
0.25 & &84.28	&62.23	 &33.16	 &18.28	 &11.21 \\
0.50 & &62.11	&32.50	 &19.84	 &13.03	 &~6.87 \\
\textbf{0.75} & &51.43	&28.95	 &19.50  &11.24  &\textbf{~6.67} \\
1~~~ & &55.84	&33.06	 &19.90	 &12.55	 &~7.74 \\
\bottomrule
\end{tabular}}
\caption{Percentage of restarted iterations in the noiseless case (Left: NLCG; Right: LBFGS).}
\label{tab:rstnoise0}
\end{table}

\begin{figure}
\centering
\subfigure[$p=0.75$ (the curves for $\kappa \ge 10^4$ overlap with the 
NLCG curve). \label{fig:sensnonoiseNLCG:fixp}]{
\includegraphics[width=0.46\textwidth]{
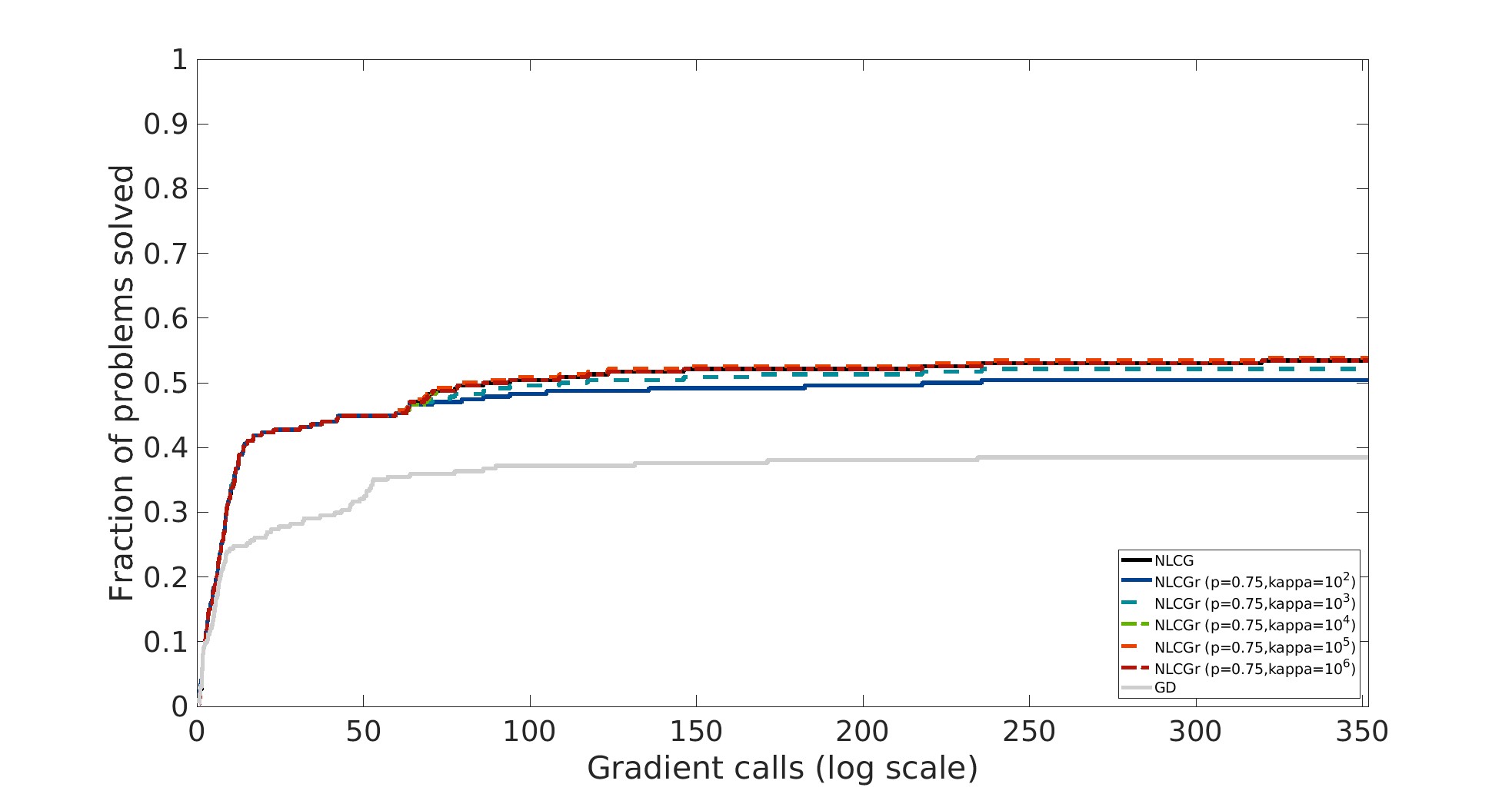}
}
\quad
\subfigure[$\kappa=10^6$ (the curves for $p \ge 0.25$ overlap with the 
NLCG curve).\label{fig:sensnonoiseNLCG:fixkappa}]{
\includegraphics[width=0.46\textwidth]{
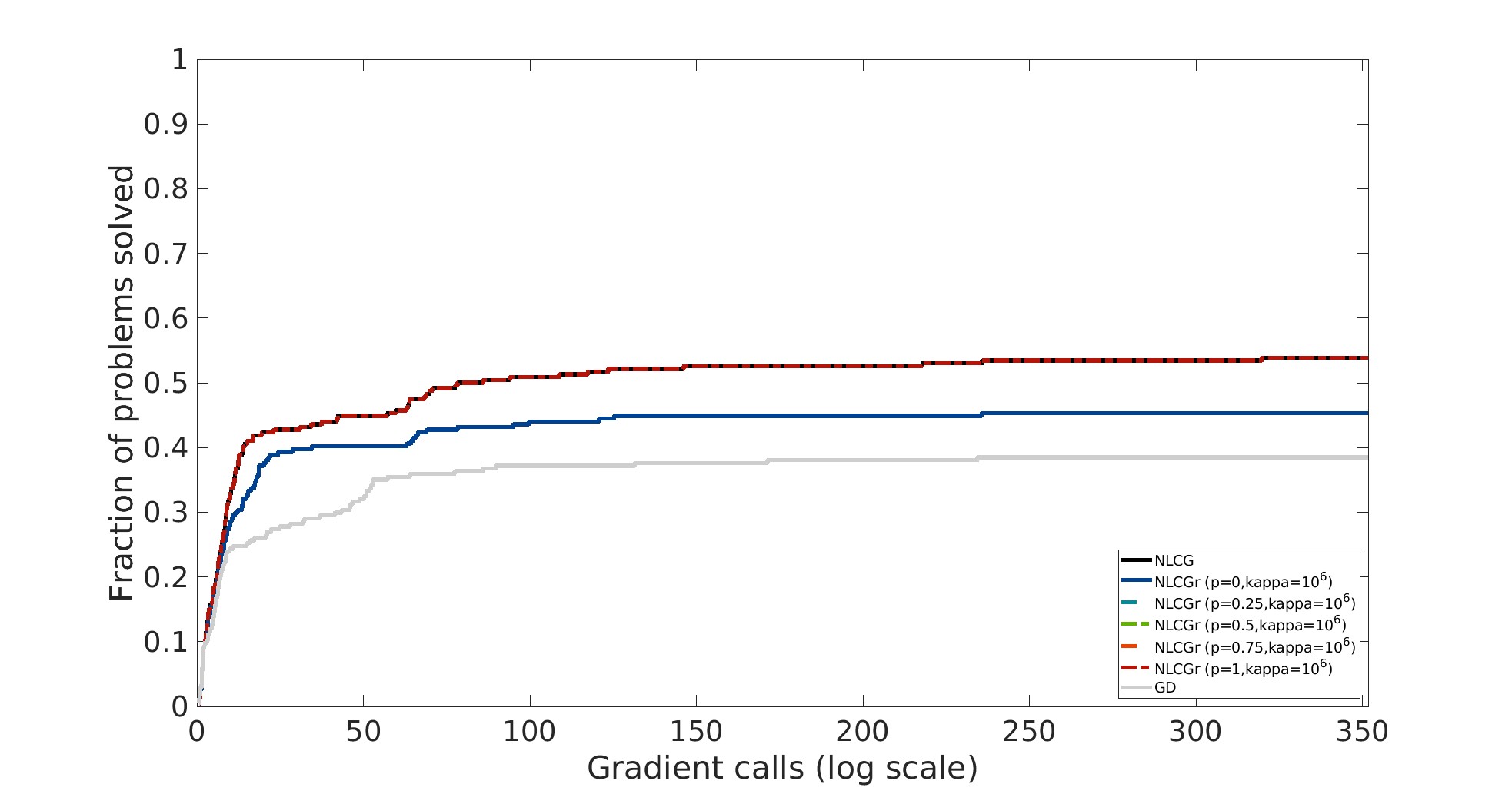}
}
\caption{Sensitivity of the performance of restarted nonlinear conjugate 
gradient (NLCGr) in absence of noise (Data profiles).}
\label{fig:sensnonoiseNLCG}
\end{figure}

\begin{figure}
\centering
\subfigure[$p=0.75$. \label{fig:sensnonoiseLBFGS:fixp}]{
\includegraphics[width=0.46\textwidth]{
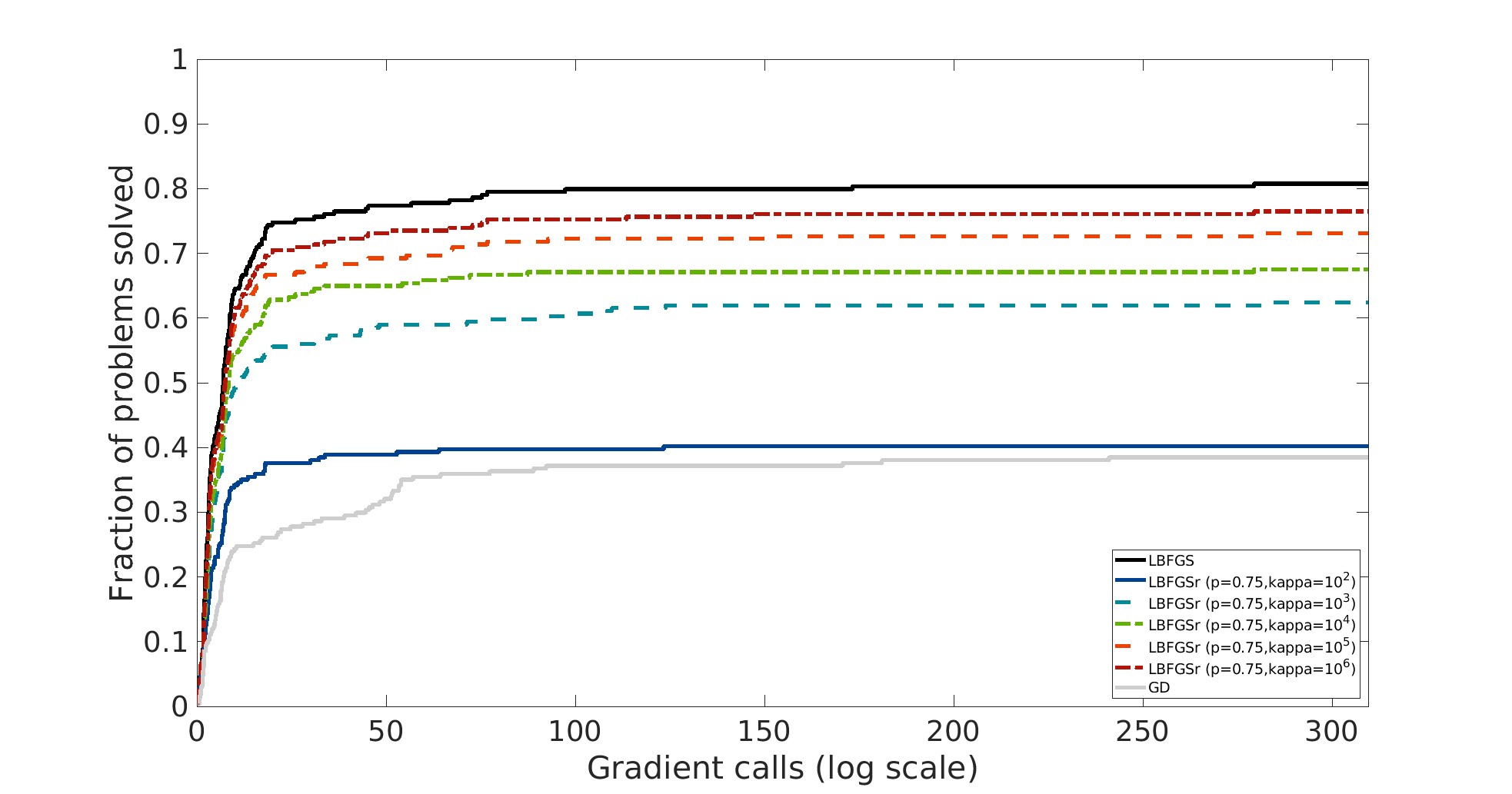}
}
\quad
\subfigure[$\kappa=10^6$.\label{fig:sensnonoiseLBFGS:fixkappa}]{
\includegraphics[width=0.46\textwidth]{
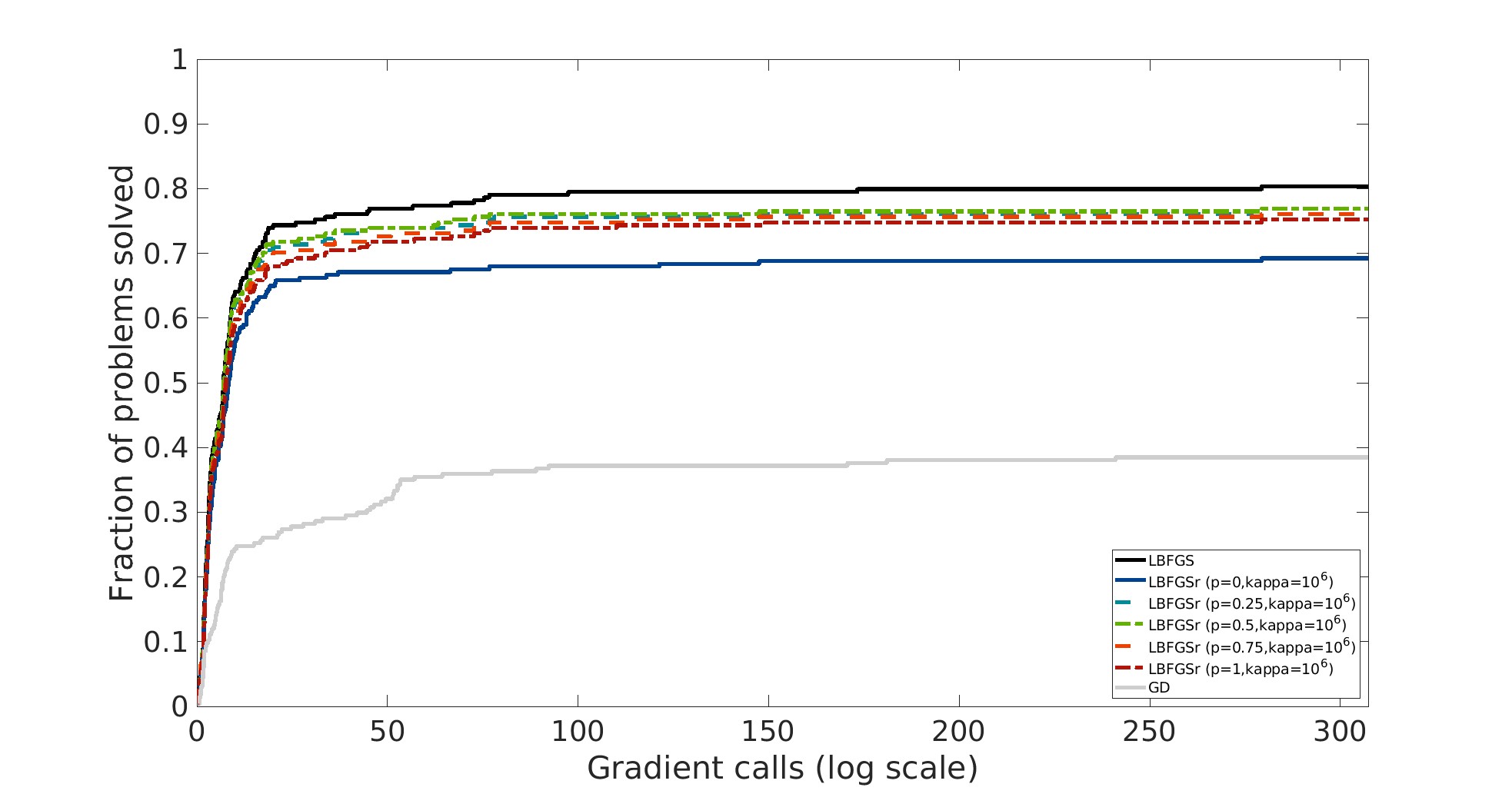}
}
\caption{Sensitivity of the performance of restarted LBFGS (LBFGSr) in absence of noise (Data profiles).}
\label{fig:sensnonoiseLBFGS}
\end{figure}

Table~\ref{tab:rstnoise0} illustrates that higher values of $\kappad$ (and 
thus smaller values of $\sigmad$) lead to fewer restarts, i.e., to algorithmic 
behavior closer to that of the non-restarted variants. A similar observation 
applies when increasing the value of $p$. Note, however, that this trend is 
more pronounced for the NLCG variants than for the LBFGS variants. To further 
highlight the impact of the restarting condition, we compare restarted variants 
using $p=0.75$ (with variable $\kappad$) and $\kappad=10^6$ 
(with variable $p$). Figure~\ref{fig:sensnonoiseNLCG} shows that the 
performance of restarted NLCG is quite close to that of standard NLCG 
except for the lowest values of $p$ and $\kappad$. On the other hand, as 
shown by Figure~\ref{fig:sensnonoiseLBFGS}, the choice of $p$ and $\kappad$ 
has a significant impact on the performance of LBFGSr, with the lowest values 
leading to a significant number of restarts and thus, performance close 
to that of gradient descent. Overall, these figures suggest that limiting the 
number of restarts leads to improved performance, closer to that of the original 
method.

\begin{figure}[]
\centering
\subfigure[Performance profile. \label{fig:compnonoise:perf}]{
\includegraphics[width=0.46\textwidth]{
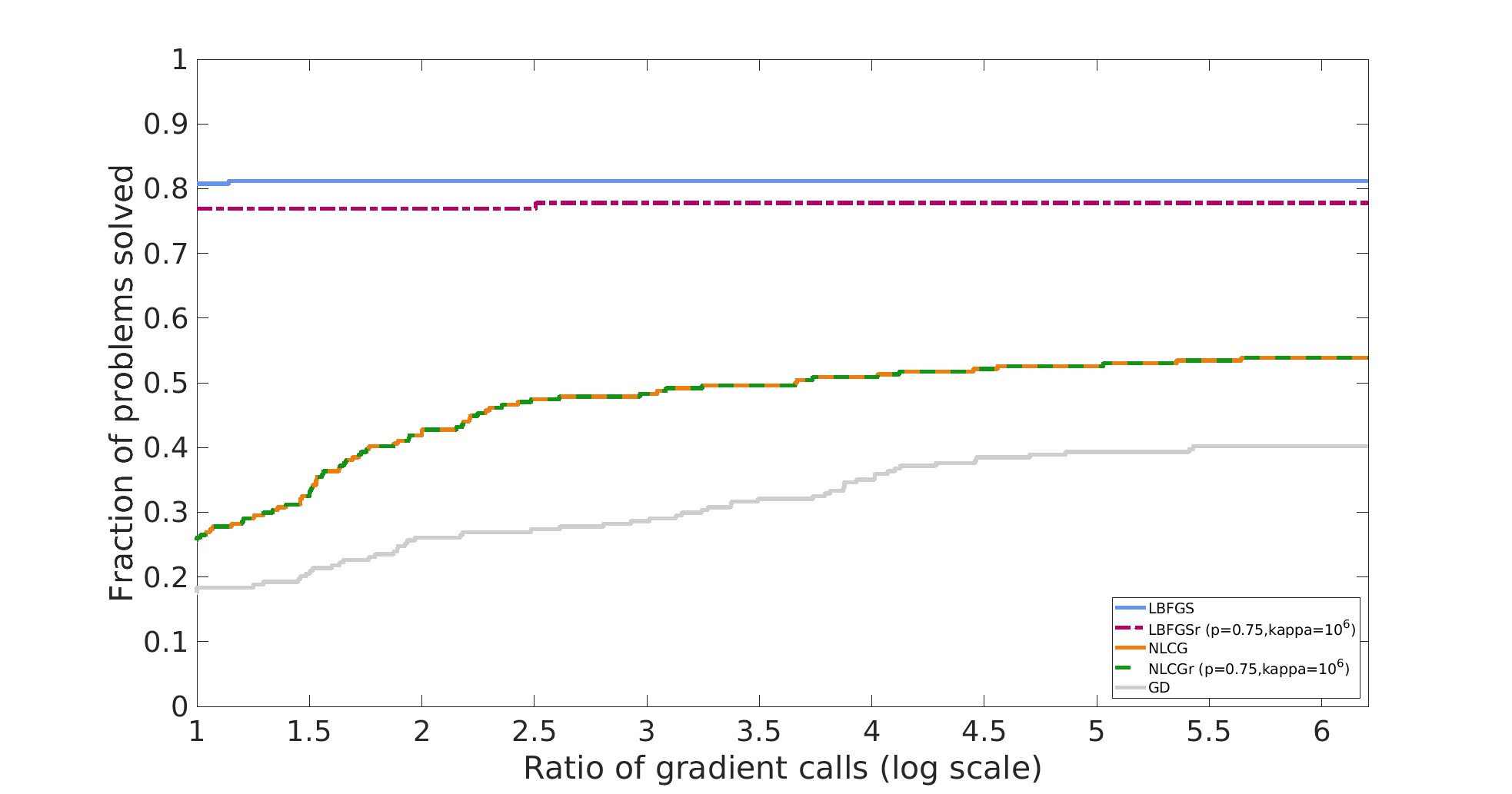}
}
\quad
\subfigure[Data profile.\label{fig:compnonoise:data}]{
\includegraphics[width=0.46\textwidth]{
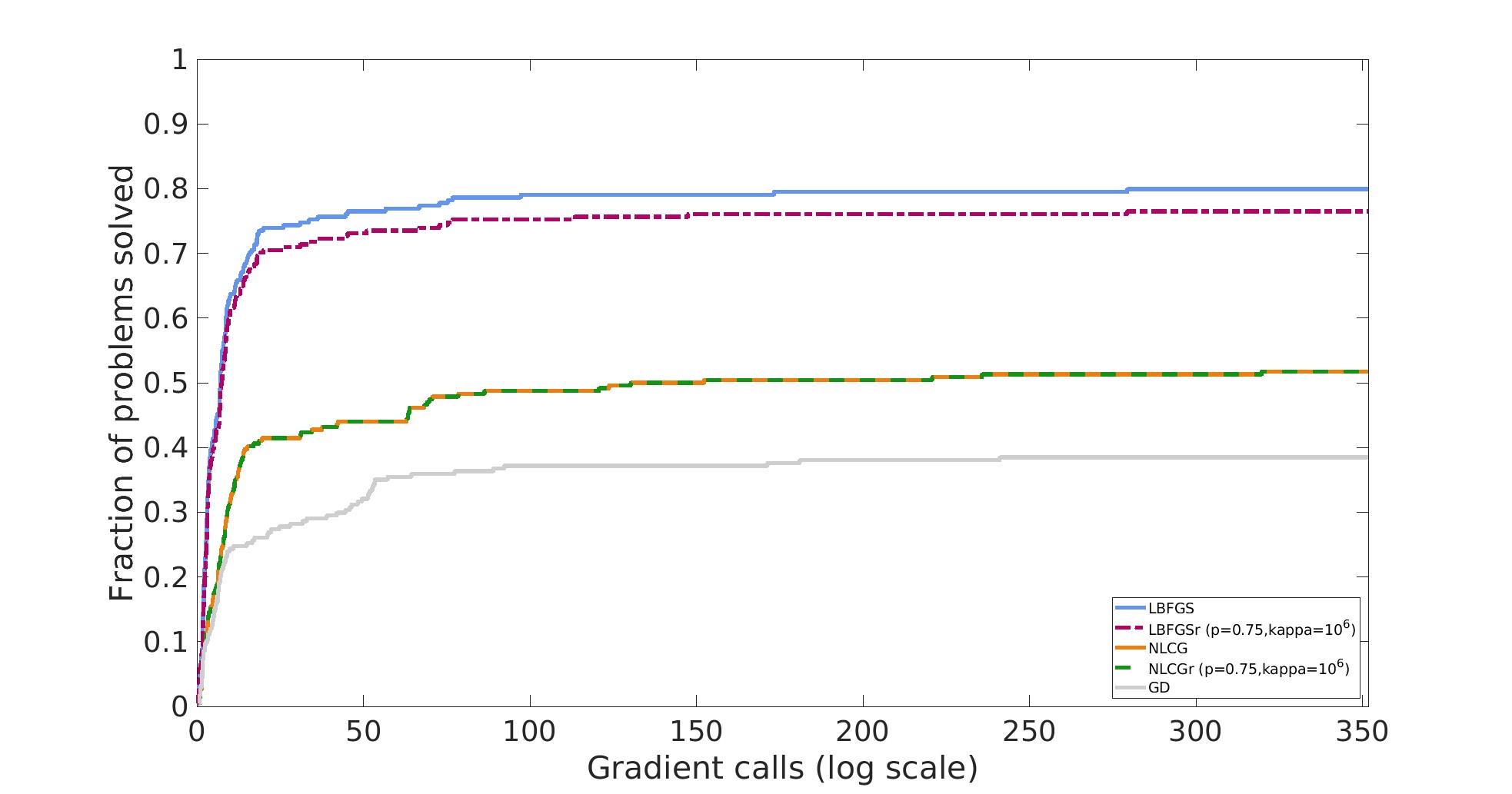}
}
\caption{Non-restarted vs. best restarted methods in the absence of noise.}
\label{fig:compnonoise}
\end{figure}

We thus compare GD, classical NLCG and LBFGS with NLCGr and LBFGSr 
using $(p=0.75,\kappa=10^{6})$. Figure~\ref{fig:compnonoise} shows that 
the two NLCG variants perform identically, whereas the performance of 
LBFGSr is slightly below that of LBFGS. These results confirm the findings 
of Chan--Renous-Legoubin and Royer in the case of 
NLCG~\cite{RChanRenousLegoubin_CWRoyer_2022}, while shedding a new light on 
the use of the restarting condition for LBFGS in the absence of noise.

\subsection{Results in noisy settings}
\label{ssec:numnoise}

In the noisy setting, our experiments reveal a higher variability in the number 
of restarts for both NLCGr and LBFGSr. Due to the randomized nature 
of the problems, we now define the best values for $p$ (resp., $\kappad$) 
according to the best performance across all possible choices for $\kappad$ 
(resp., $p$). Note that our stopping criterion is less stringent than in the 
noiseless setting.

Detailed statistics regarding restarts are given in 
Appendix~\ref{app:restarttables}, and Table~\ref{tab:rstnoiseall} 
summarizes the best configurations in terms of $p$ and $\kappad$ (recall that 
we set $\sigmad=\kappad^{-1}$).
For LBFGSr, using a large value of $\kappad$ is critical, which we attribute 
to the directions being of relatively high norm compared to that of the 
gradient. Lower values of $p$ appear preferable as the noise increases. 
For NLCGr, both the best value of $\kappad$ and that of $p$ decrease as the 
noise increases. Note these observations are consistent with our  analysis, 
in that the selected values lead to better complexity bounds.

\begin{table}[]
\centering
\begin{tabular}{lrrrr}
\toprule
Method &$\epsf=10^{-8}$ &$\epsf=10^{-4}$ &$\epsf=10^{-2}$ 
&$\epsf=10^{-1}$ \\
\midrule
NLCGr       &$(~0.75,10^{6})$  &$(0.75,10^{5})$ &$(0,10^{3})$ 
&$(0.25,10^{3})$ \\
LBFGSr     &$(0.75,10^{6})$ &$(~~1,10^{6})$ &$(0.5,10^{6})$ 
&$(~~0,10^{6})$ \\
\bottomrule
\end{tabular}
\caption{Values of $p$ and $\kappad$ with the least percentage of restarted 
iterations.}
\label{tab:rstnoiseall}
\end{table}

\begin{figure}
\centering
\subfigure[$p=0.75$ (the curves for $\kappa \ge 10^3$ overlap with the 
NLCG curve). \label{fig:sensnoise8NLCG:fixp}]{
\includegraphics[width=0.46\textwidth]{
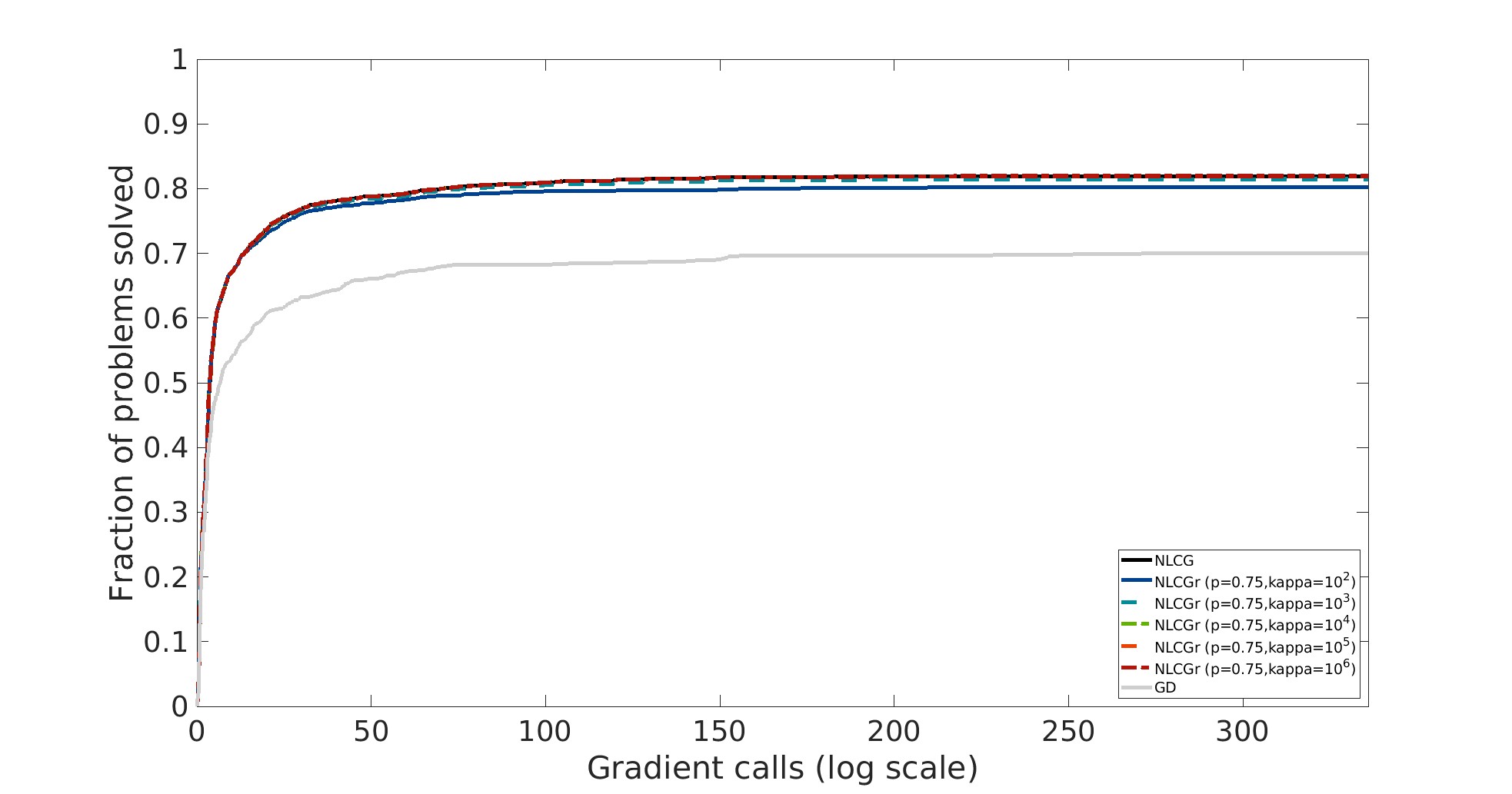}
}
\quad
\subfigure[$\kappa=10^6$ (all NLCG curves overlap with the 
NLCG curve).\label{fig:sensnoise8NLCG:fixkappa}]{
\includegraphics[width=0.46\textwidth]{
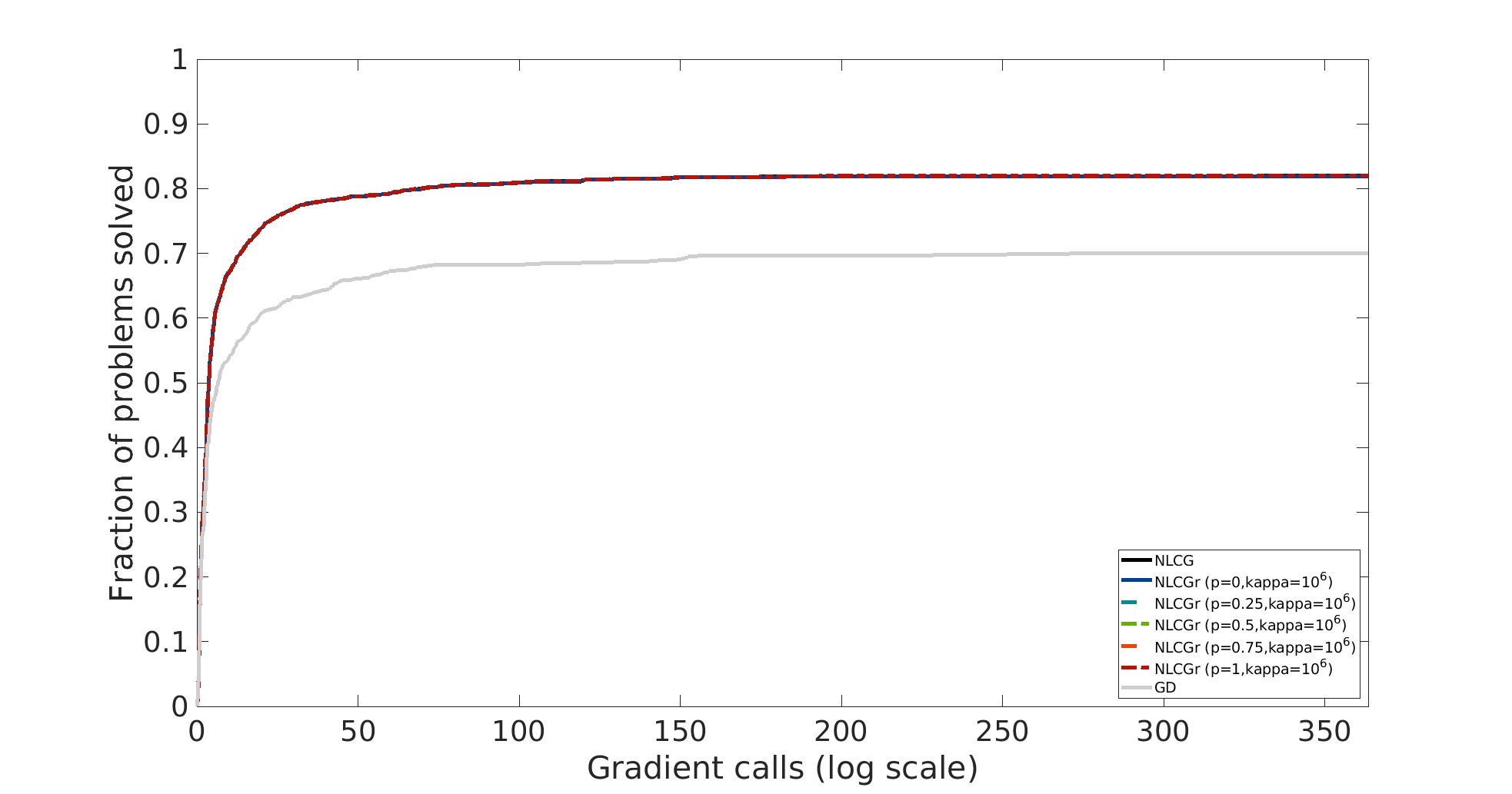}
}
\caption{Sensitivity of the performance of restarted nonlinear conjugate 
gradient (NLCGr) with $\epsf=10^{-8}$ (Data profiles).}
\label{fig:sensnoise8NLCG}
\end{figure}

\begin{figure}
\centering
\subfigure[$p=0.75$. \label{fig:sensnoise8LBFGS:fixp}]{
\includegraphics[width=0.46\textwidth]{
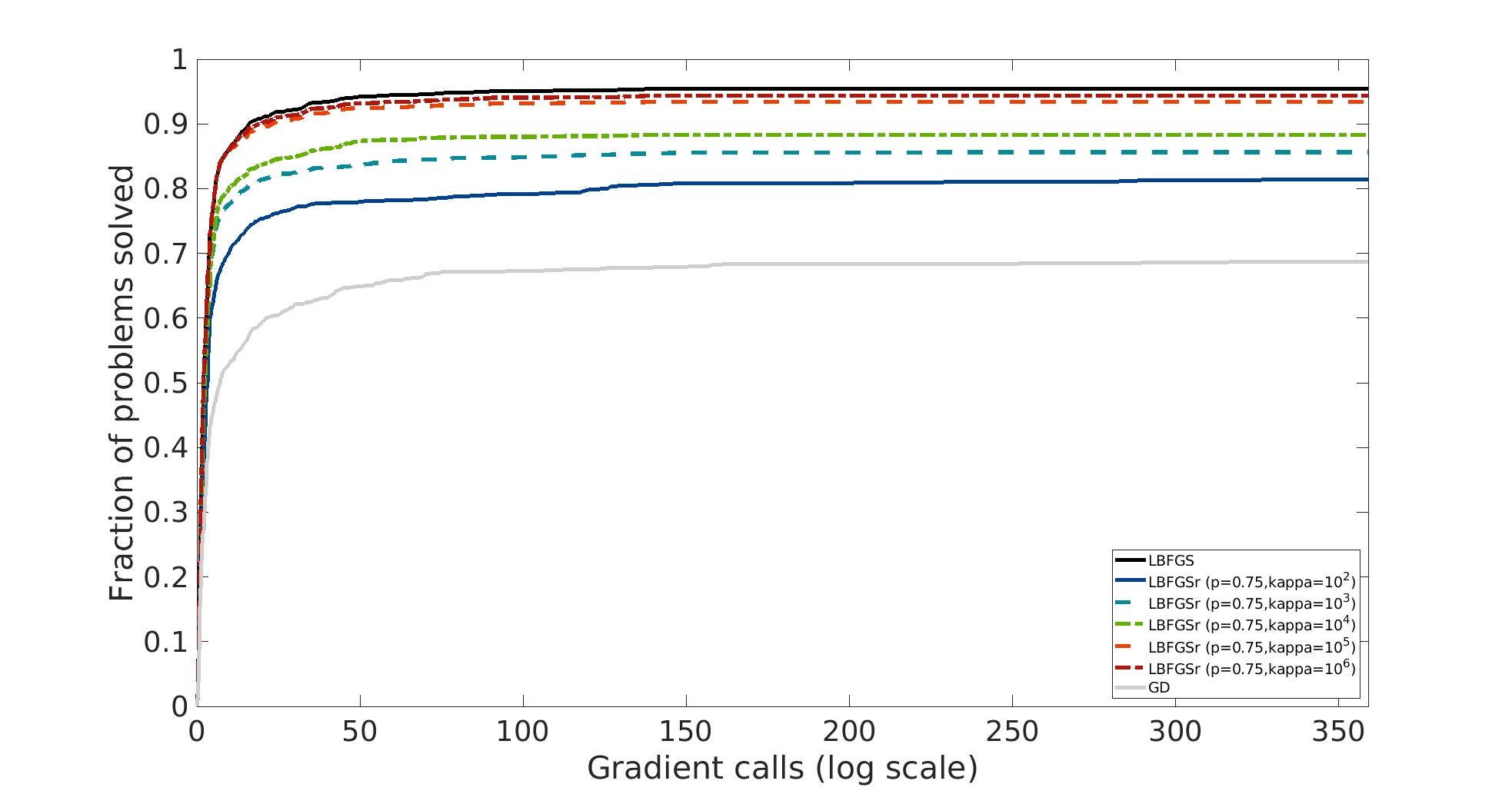}
}
\quad
\subfigure[$\kappa=10^6$ (the curves $p \in \{0.25,0.5\}$ and 
$p \in \{0.75,1\}$ respectively overlap).\label{fig:sensnoise8LBFGS:fixkappa}]{
\includegraphics[width=0.46\textwidth]{
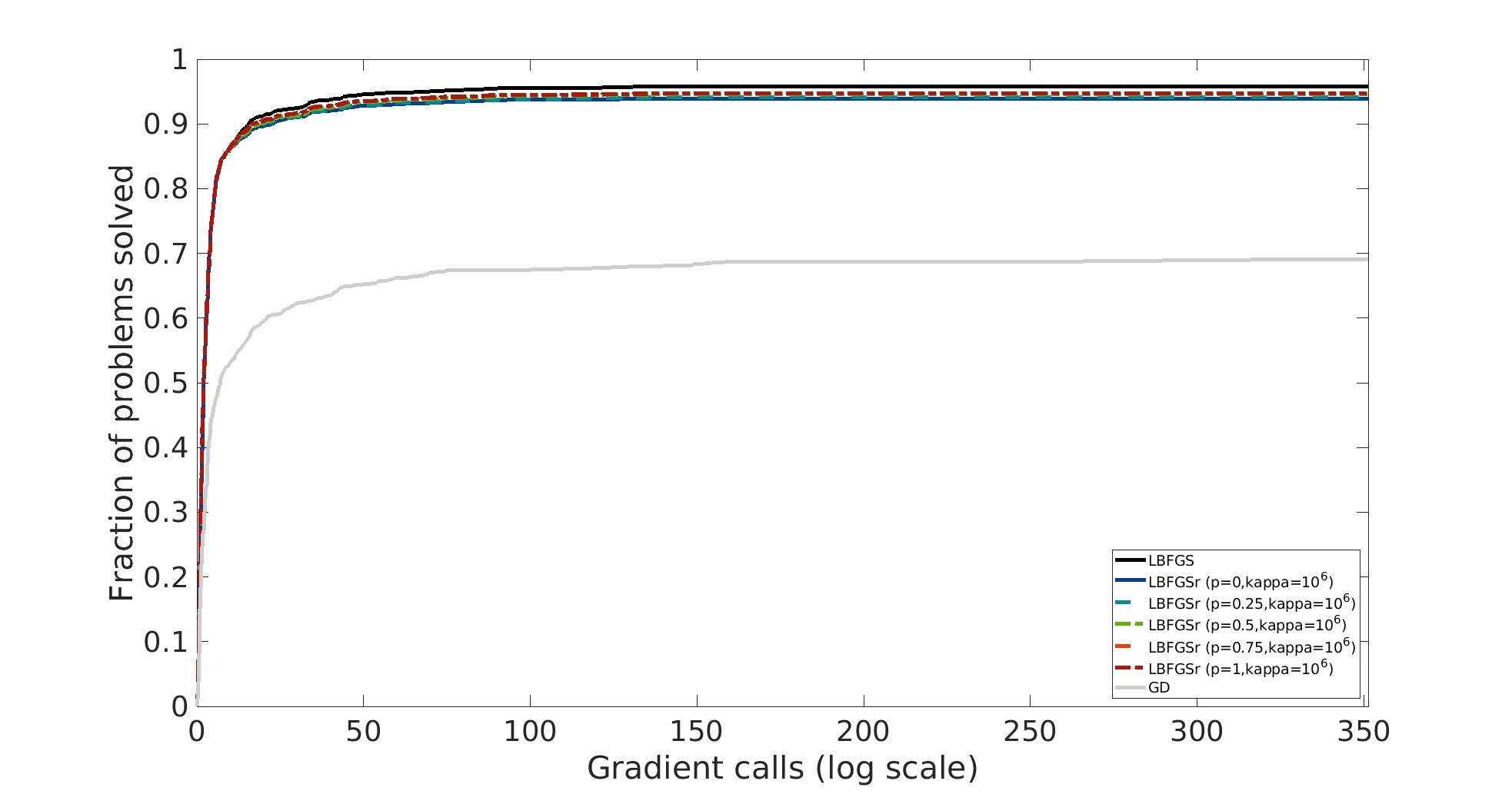}
}
\caption{Sensitivity of the performance of restarted LBFGS (LBFGSr) with $\epsf=10^{-8}$ (Data profiles).}
\label{fig:sensnoise8LBFGS}
\end{figure}

To further motivate these parameter choices, we conducted a sensitivity 
analysis for the lowest and highest noise levels. When the noise level is 
relatively low, Figures~\ref{fig:sensnoise8NLCG} 
and~\ref{fig:sensnoise8LBFGS} yield the same conclusions as in the noiseless 
case. Indeed, the NLCGr variants coincide with NLCG for sufficiently high values of 
the parameters. A similar conclusion holds for LBFGSr, although we notice that 
all variants are much closer than in the noiseless case, which we attribute to 
the lower number of overall restarts performed by those variants (see 
Table~\ref{tab:rstnoise1em8}, Appendix \ref{app:restarttables} for detailed 
statistics).


\begin{figure}
\centering
\subfigure[$p=0.25$. \label{fig:sensnoise1NLCG:fixp}]{
\includegraphics[width=0.46\textwidth]{
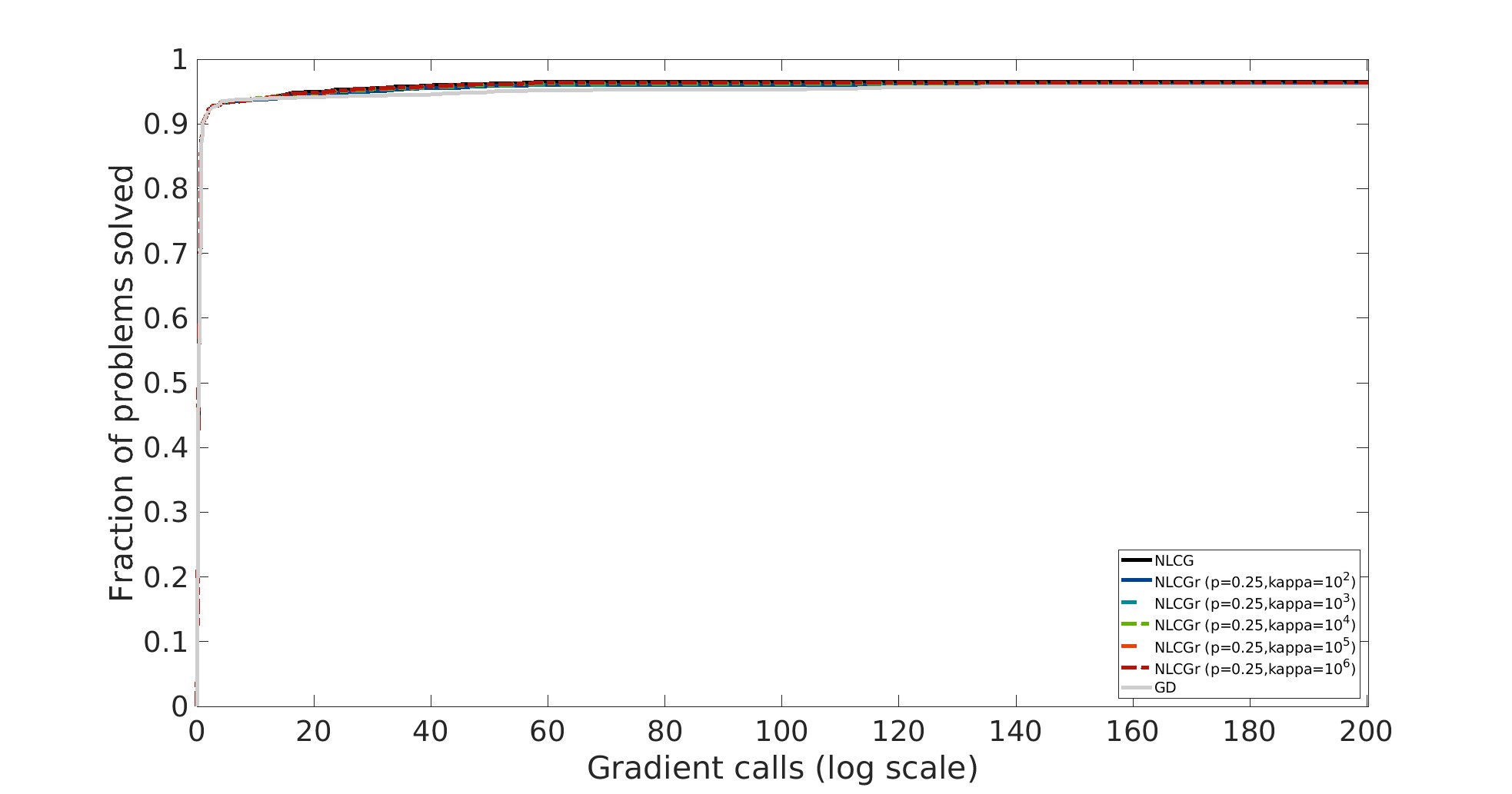}
}
\subfigure[$\kappa=10^3$.\label{fig:sensnoise1NLCG:fixkappa}]{
\includegraphics[width=0.46\textwidth]{
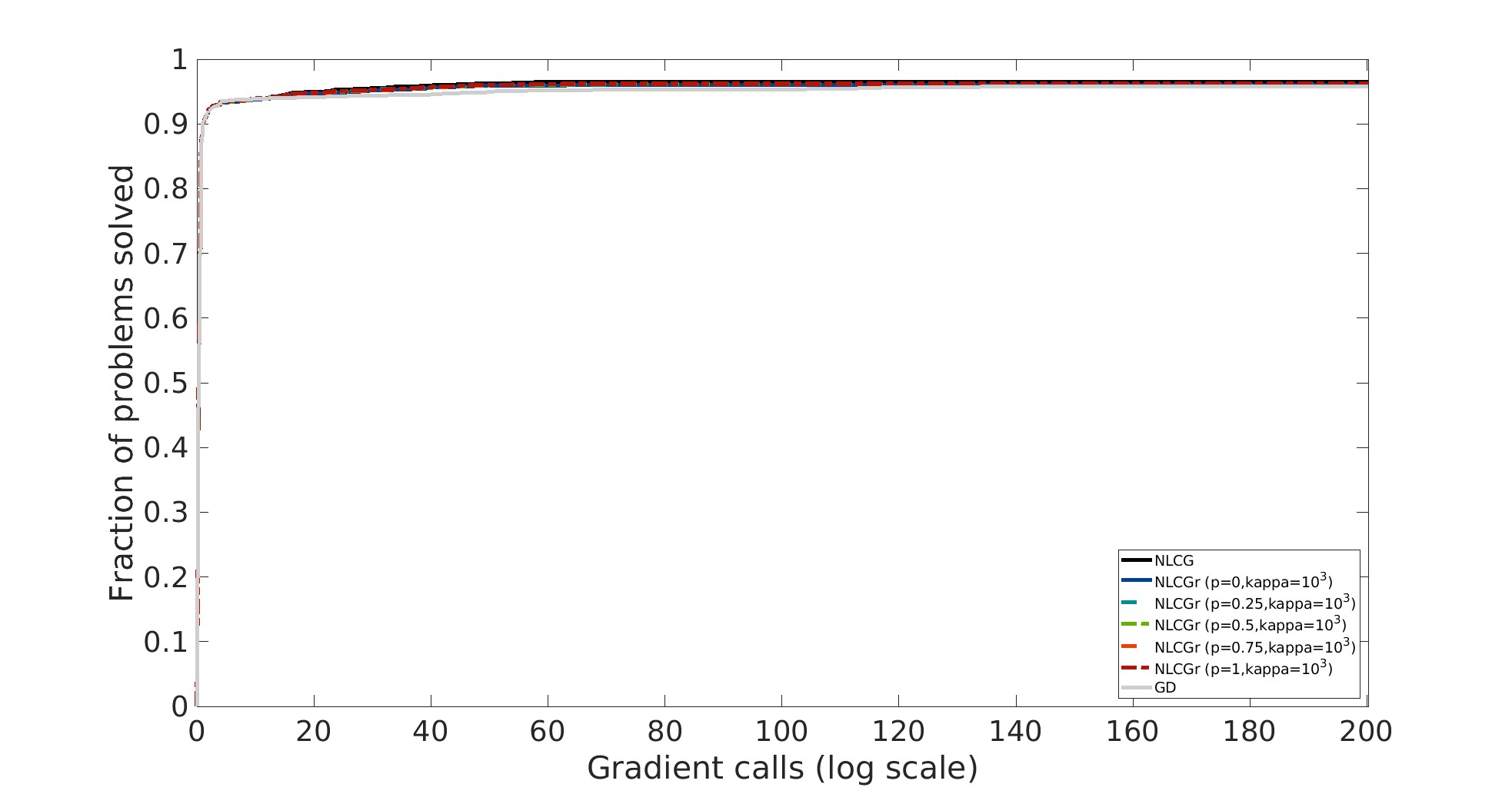}
}
\quad
\caption{Sensitivity of the performance of restarted nonlinear conjugate 
gradient with $\epsf=10^{-1}$ (Data profiles).}
\label{fig:sensnoise1NLCG}
\end{figure}

\begin{figure}
\centering
\subfigure[$p=0$. \label{fig:sensnoise1LBFGS:fixp}]{
\includegraphics[width=0.46\textwidth]{
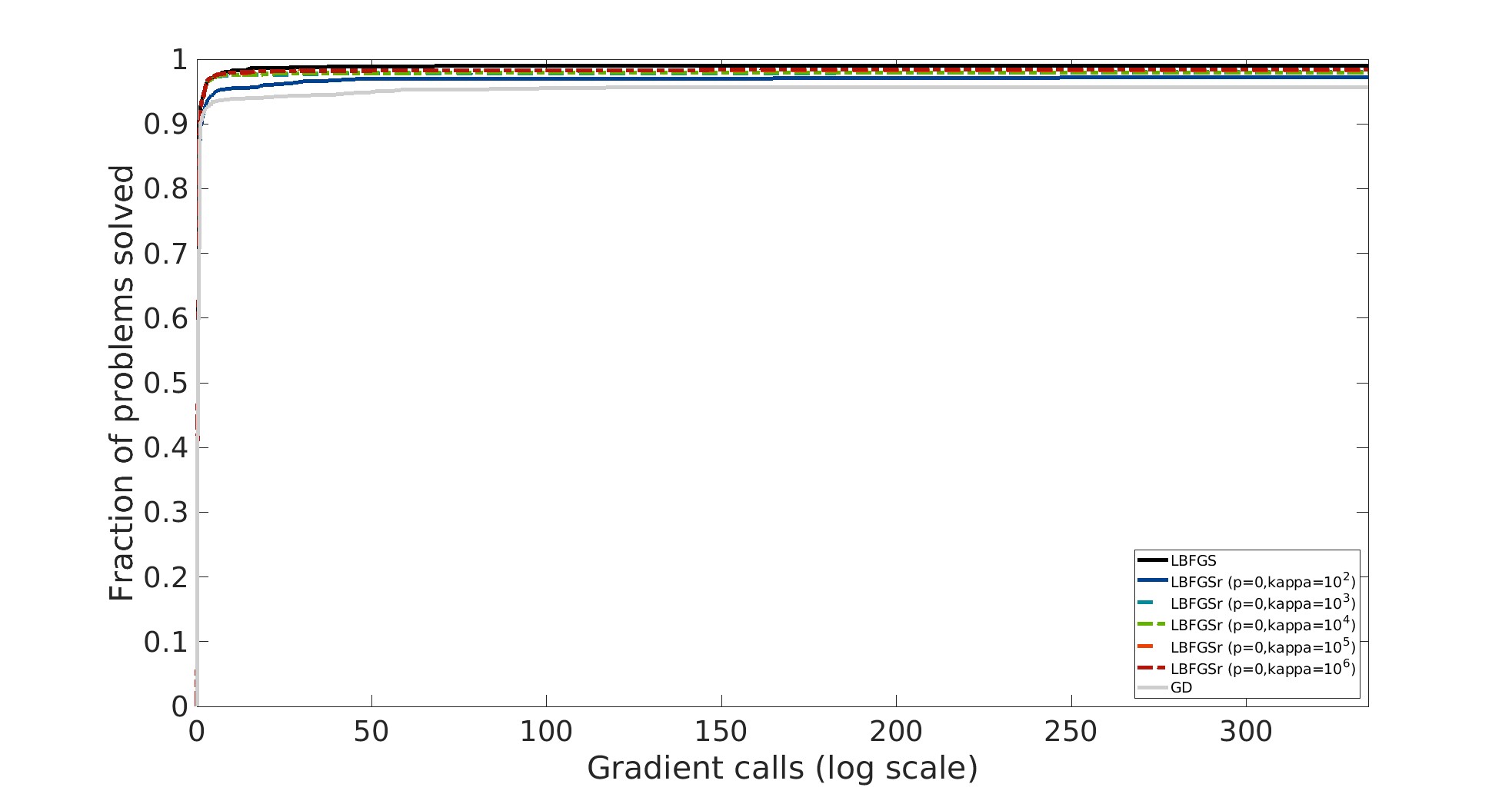}
}
\subfigure[$\kappa=10^6$.\label{fig:sensnoise1LBFGS:fixkappa}]{
\includegraphics[width=0.46\textwidth]{
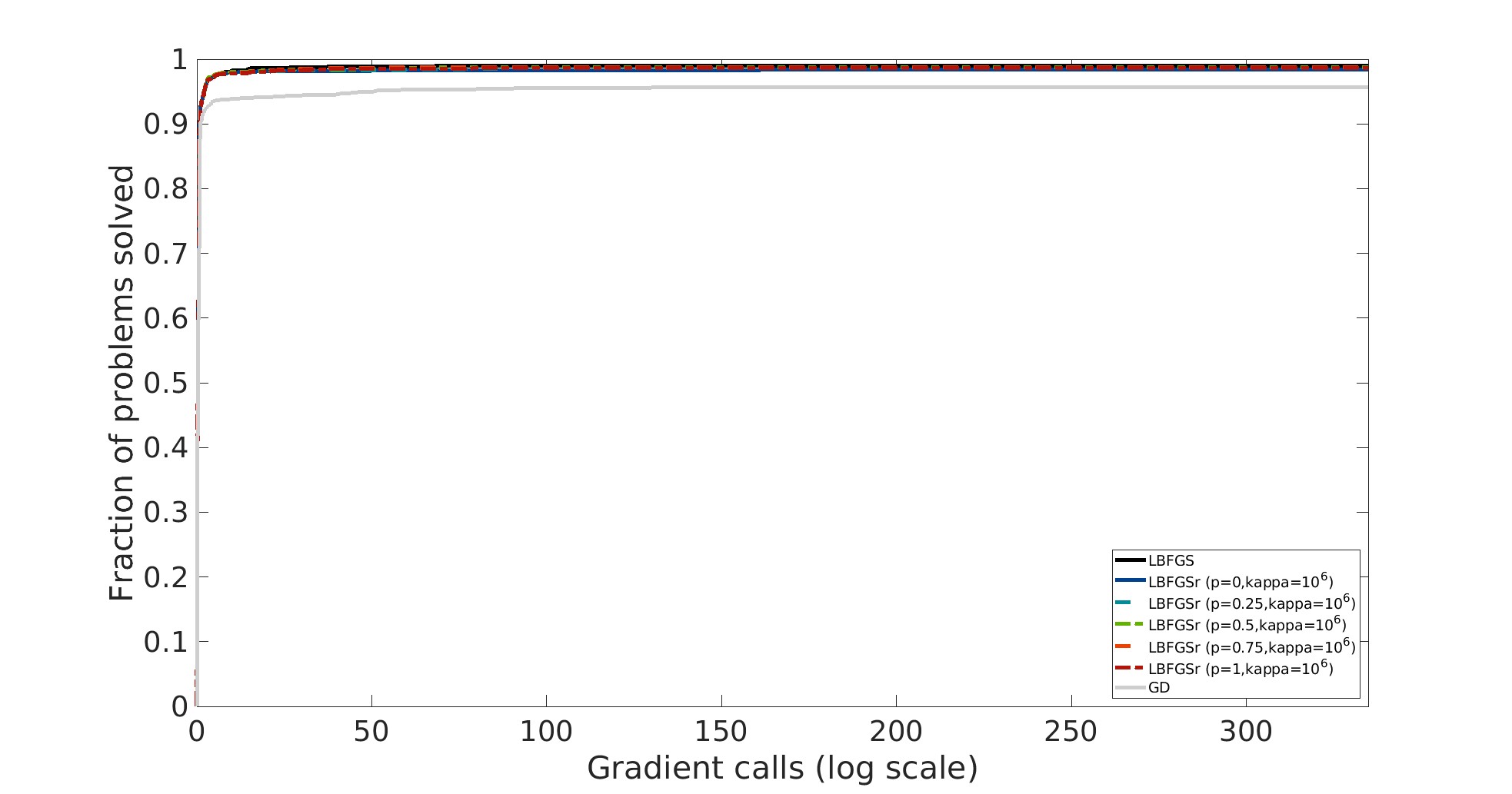}
}
\quad
\caption{Sensitivity of the performance of restarted LBFGS (LBFGSr) with $\epsf=10^{-1}$ (Data profiles).}
\label{fig:sensnoise1LBFGS}
\end{figure}

For the higher noise level (Figures~\ref{fig:sensnoise1NLCG} 
and~\ref{fig:sensnoise1LBFGS}), we notice that the curves are much closer to 
one another, for both NLCG and LBFGS. This behavior is largely explained by 
the convergence criterion used in that case. Indeed, only a low precision on 
the gradient is required, which all methods manage to satisfy rather easily. 
The number of restarts is also quite similar for all NLCGr variants, with 
more variability for LBFGSr variants (see Table~\ref{tab:rstnoise1em1}, 
Appendix~\ref{app:restarttables} for detailed statistics). Overall, our 
experiments suggest that using a restarted variant in the presence of noise 
will lead to performance comparable to that of a non-restarted variant.

\begin{figure}[h!]
\centering
\subfigure[Performance profile. \label{fig:compnoise8:perf}]{
\includegraphics[width=0.46\textwidth]{
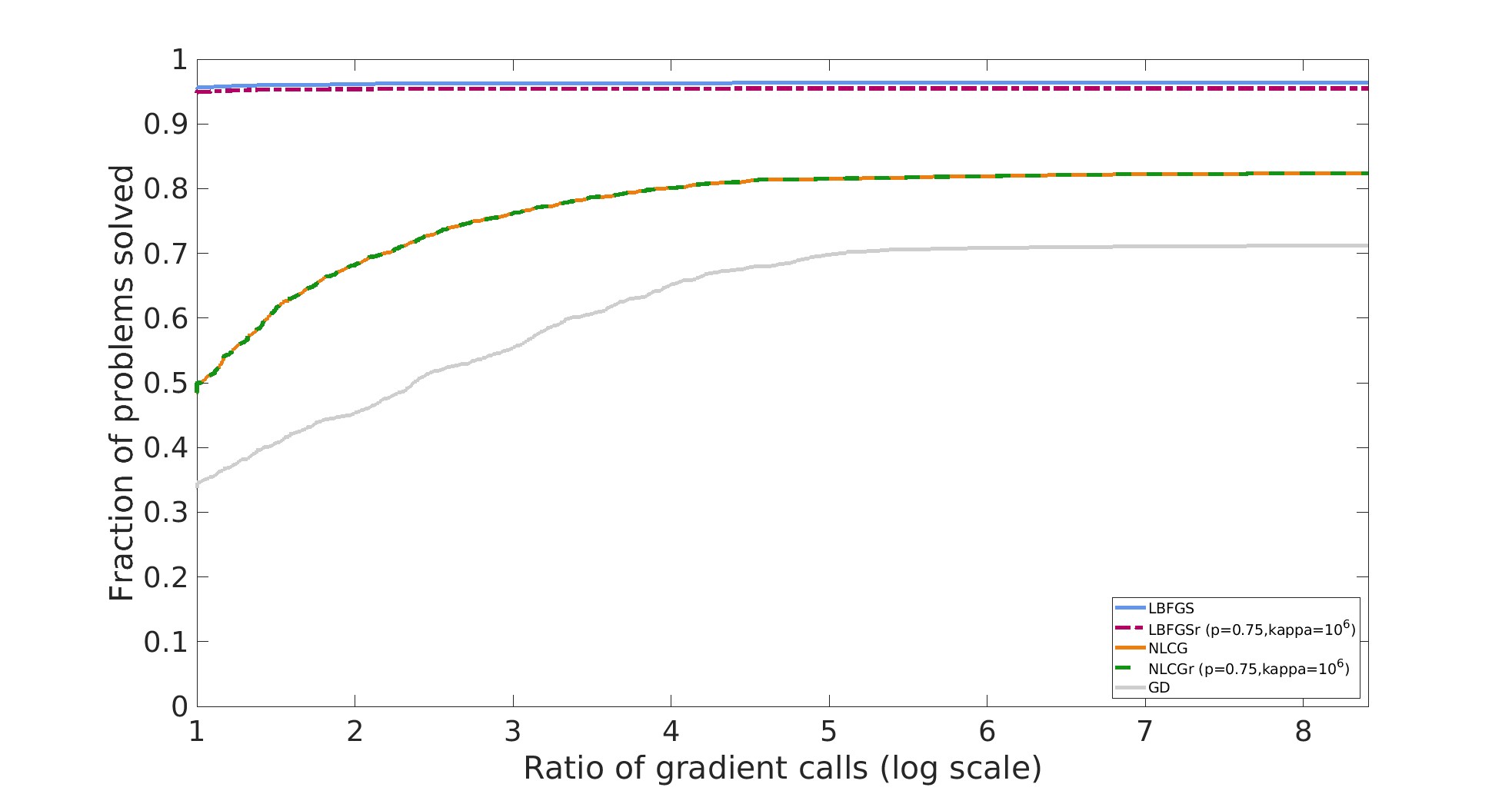}
}
\quad
\subfigure[Data profile.\label{fig:compnoise8:data}]{
\includegraphics[width=0.46\textwidth]{
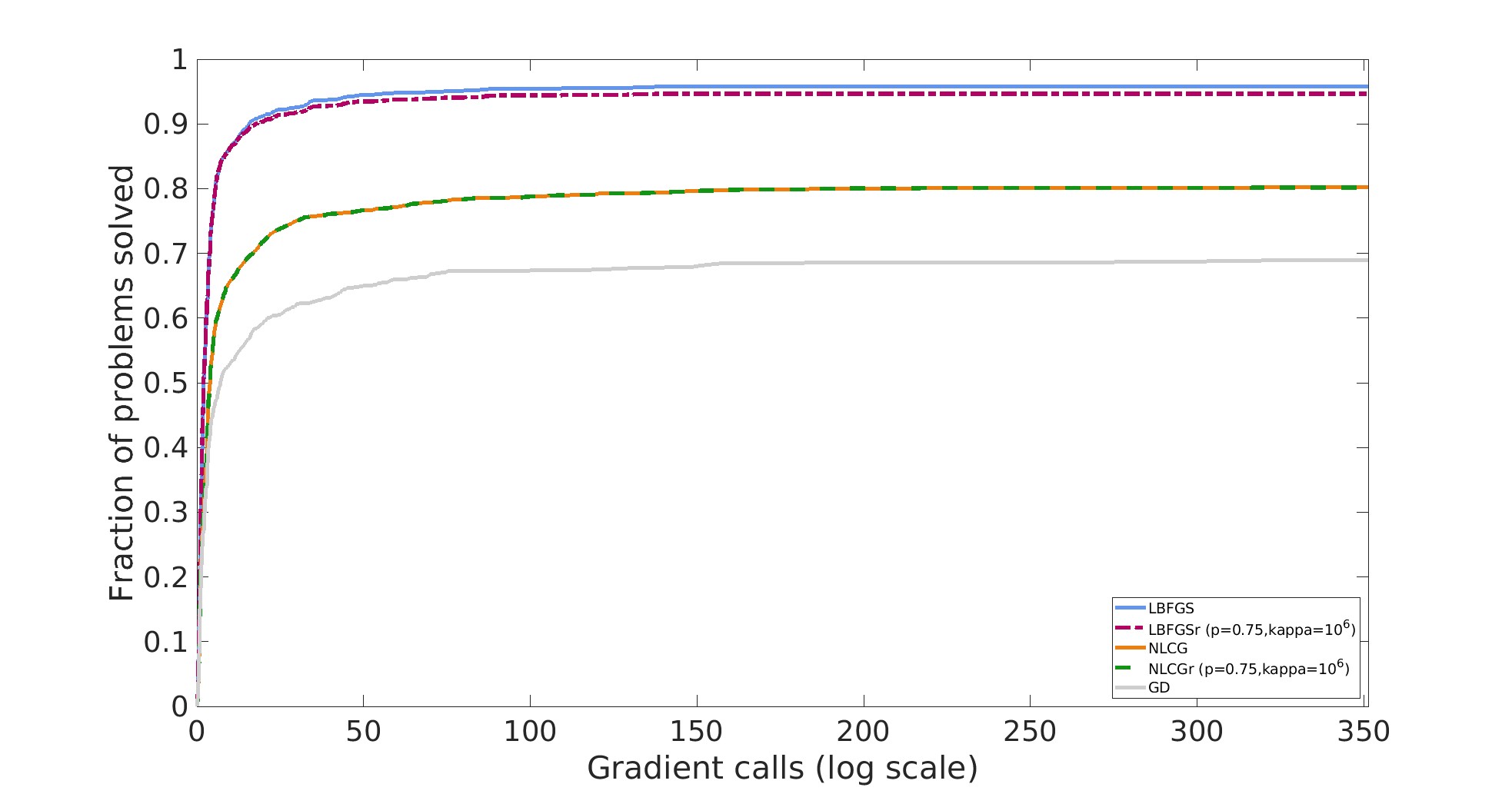}
}
\caption{Non-restarted vs. best restarted methods with $\epsf=10^{-8}$.}
\label{fig:compnoise8}
\end{figure}

\begin{figure}[h!]
\centering
\subfigure[Performance profile. \label{fig:compnoise4:perf}]{
\includegraphics[width=0.46\textwidth]{
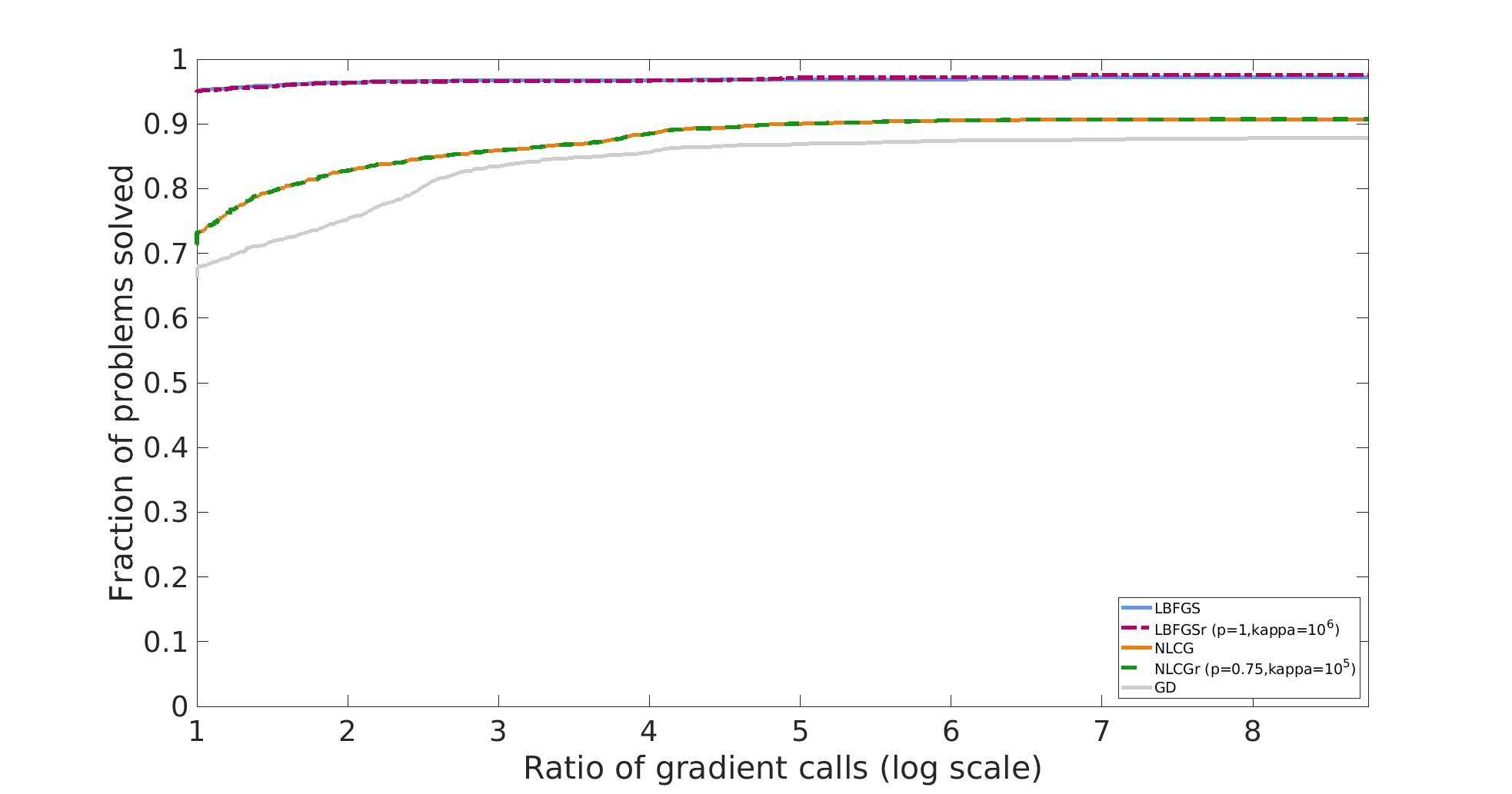}
}
\quad
\subfigure[Data profile.\label{fig:compnoise4:data}]{
\includegraphics[width=0.46\textwidth]{
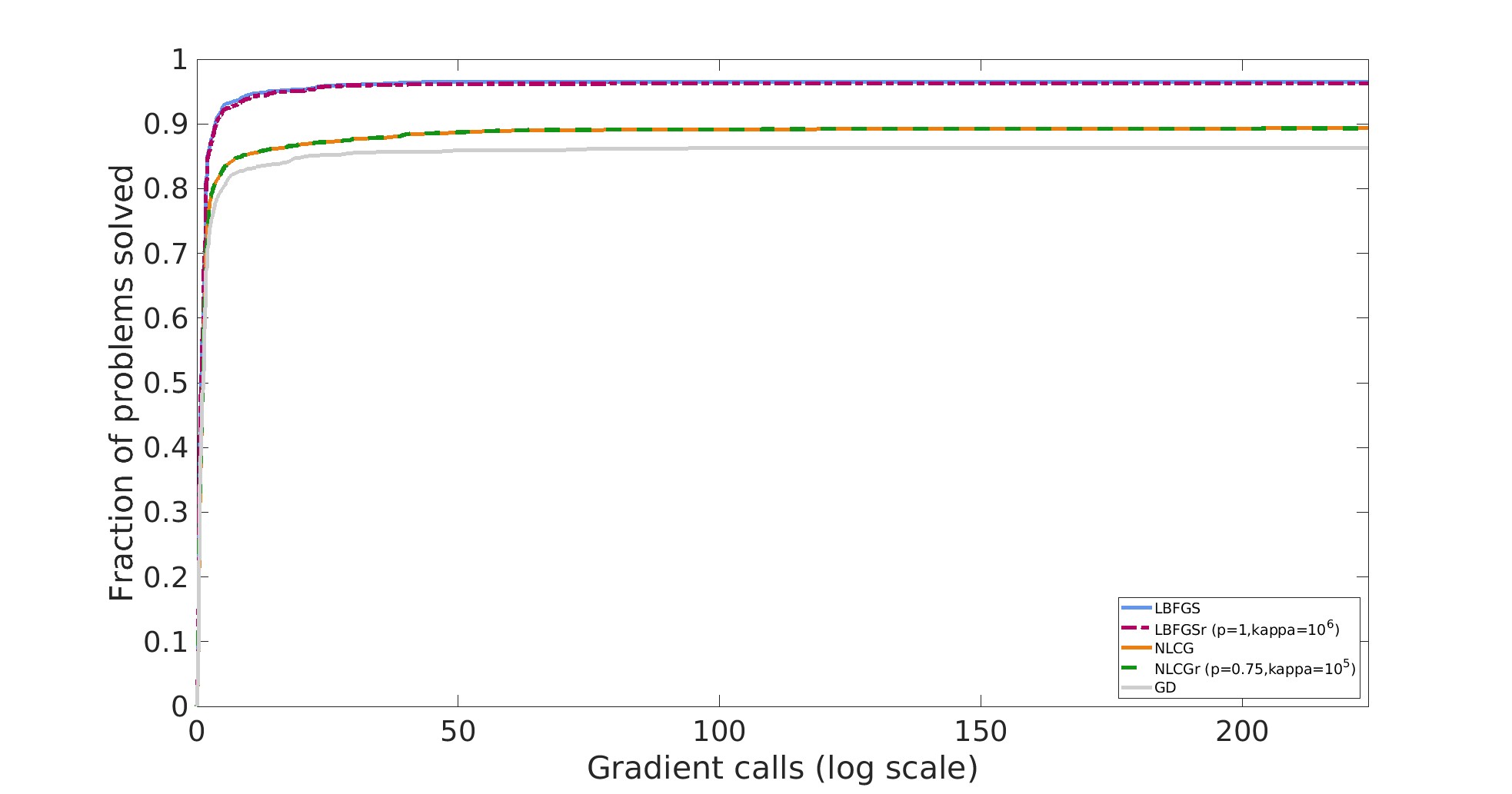}
}
\caption{Non-restarted vs.best restarted methods with $\epsf=10^{-4}$.}
\label{fig:compnoise4}
\end{figure}

\begin{figure}[h!]
\centering
\subfigure[Performance profile. \label{fig:compnoise2:perf}]{
\includegraphics[width=0.46\textwidth]{
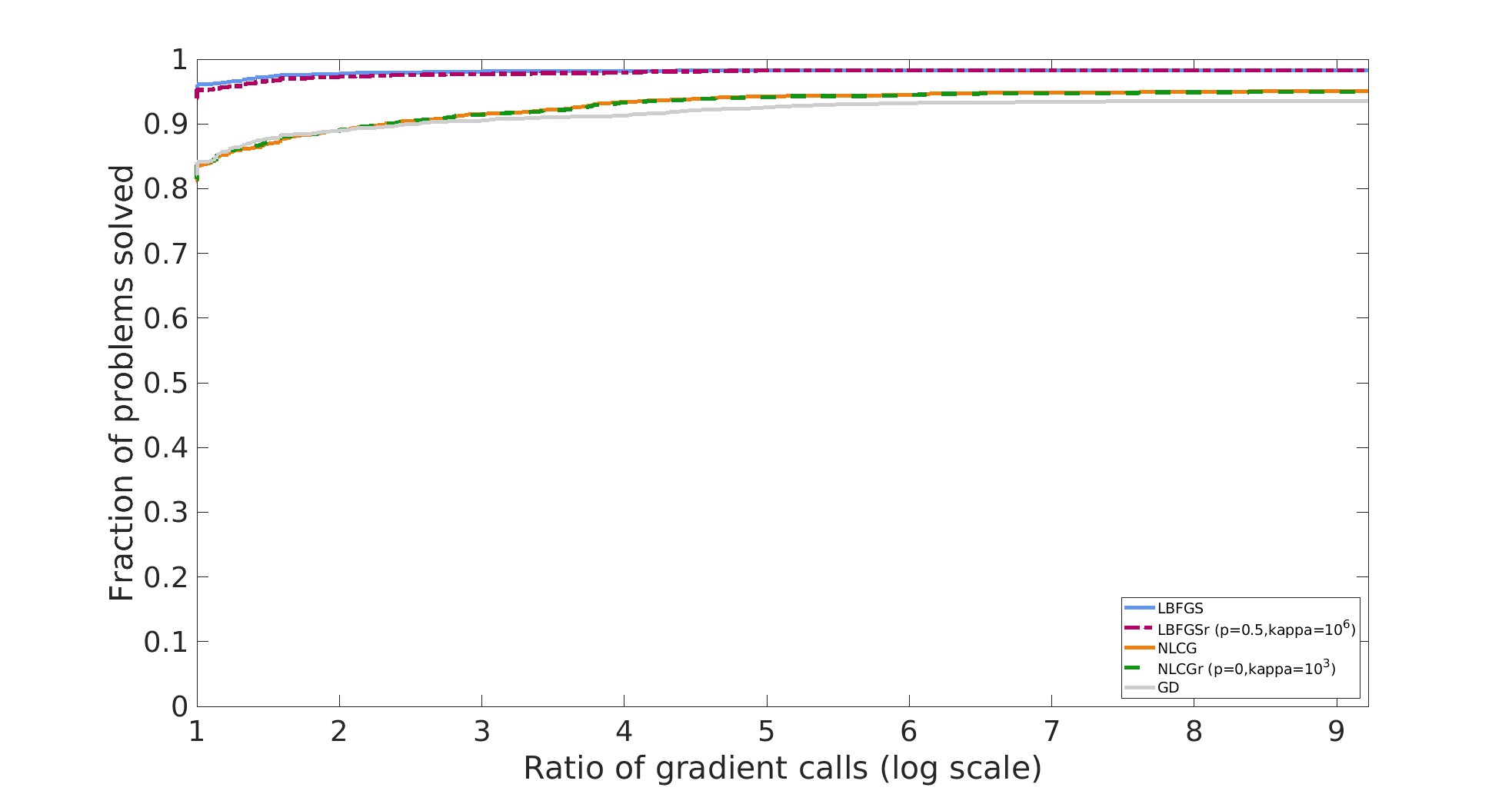}
}
\quad
\subfigure[Data profile.\label{fig:compnoise2:data}]{
\includegraphics[width=0.46\textwidth]{
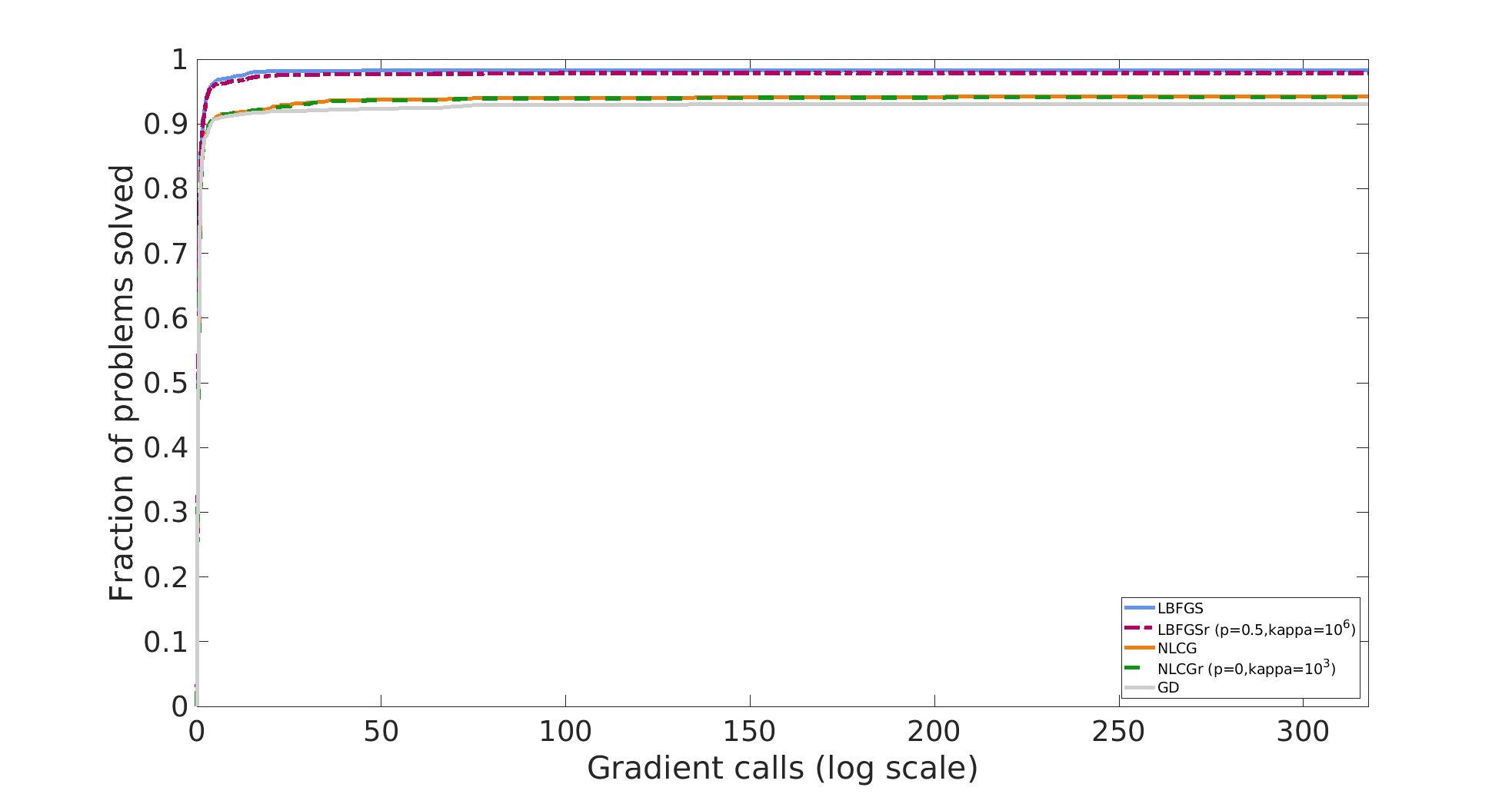}
}
\caption{Non-restarted vs.best restarted methods with $\epsf=10^{-2}$.}
\label{fig:compnoise2}
\end{figure}

\begin{figure}[h!]
\centering
\subfigure[Performance profile. \label{fig:compnoise1:perf}]{
\includegraphics[width=0.46\textwidth]{
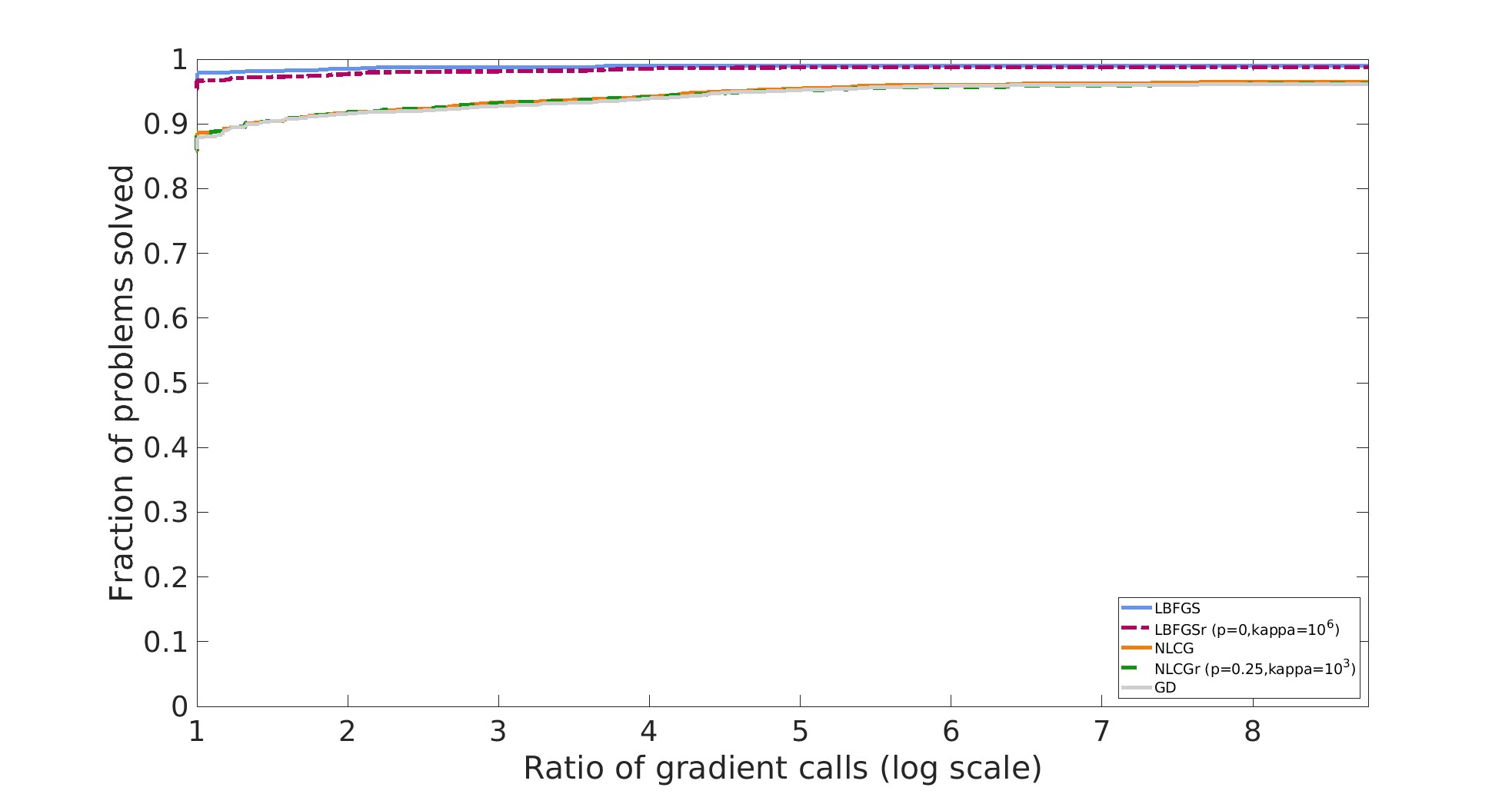}
}
\quad
\subfigure[Data profile.\label{fig:compnoise1:data}]{
\includegraphics[width=0.46\textwidth]{
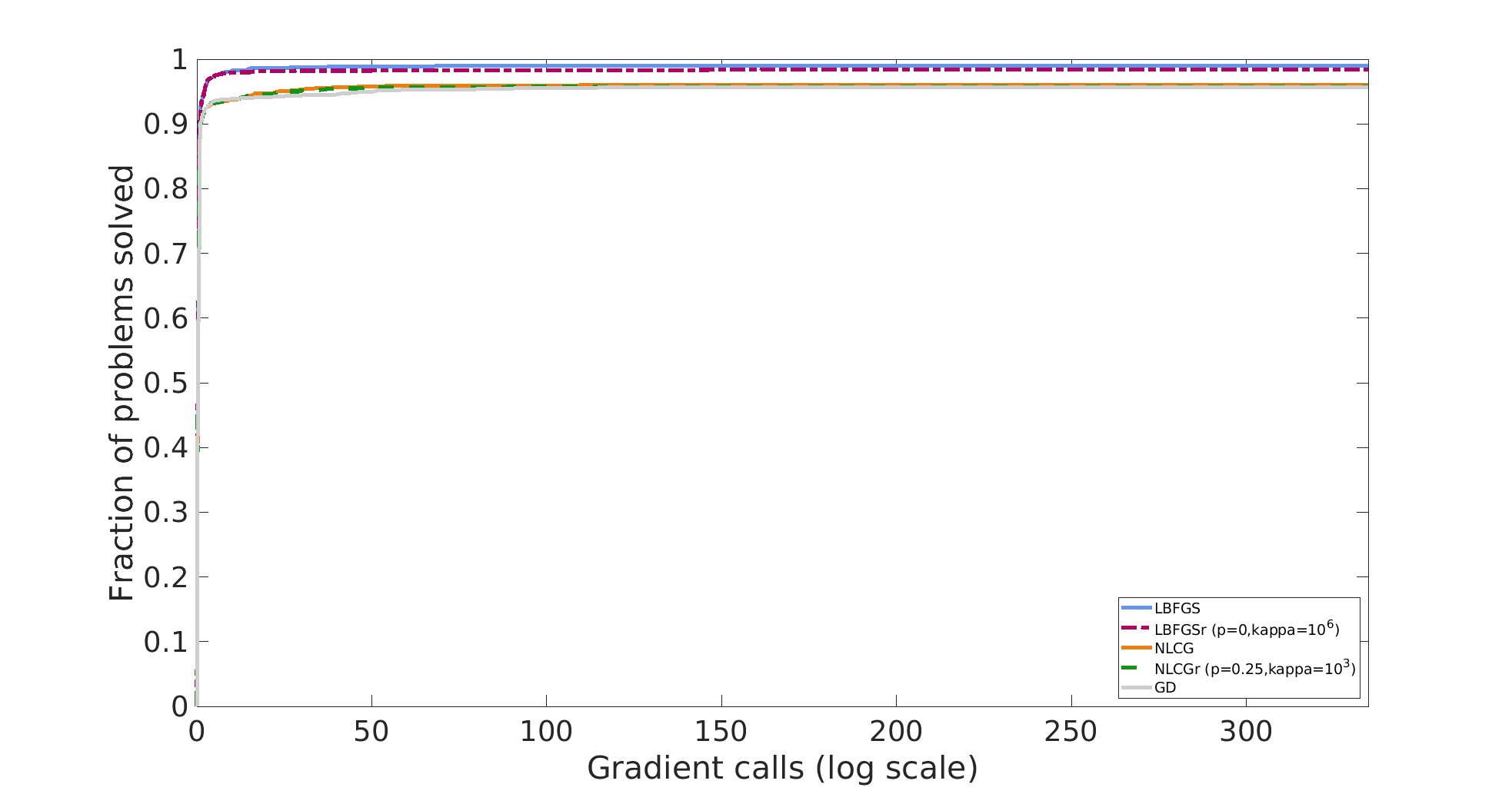}
}
\caption{Non-restarted vs.best restarted methods with $\epsf=10^{-1}$.}
\label{fig:compnoise1}
\end{figure}

To further confirm this observation, we compare GD, NLCG, LBFGS and the two 
best restarted variants according to Table~\ref{tab:rstnoiseall} for all 
noise levels. Figures~\ref{fig:compnoise8} to~\ref{fig:compnoise1} show 
performance and data profiles for four noise levels. Restarted variants are 
typically overlapping with their non-restarted counterparts, as in the 
noiseless case. As the noise increases, we notice that all curves get closer 
to one another. In particular, gradient descent matches the performance of 
NLCG variants for $\epsf \ge 10^{-2}$, and the difference between the LBFGS 
variants and others gets less pronounced as the noise increases. Still, we note 
that the LBFGS variants exhibit the best performance overall, confirming the 
interest of such approaches even in noisy settings.

\section{Conclusion}
\label{sec:conc}

We have proposed a line search framework for noisy optimization problems in 
which advanced search direction techniques revert to (inexact) gradient steps 
when a restarting condition is triggered. These conditions are instrumental to 
obtaining complexity guarantees, and our analysis highlights that the 
percentage of restarted iterations has an impact on complexity bounds. The 
proportion of restarted iterations is also a major factor in the numerical 
performance. In particular, the performance of LBFGS variants is highly 
sensitive to the restart condition parameters. When calibrated appropriately, 
however, a restarted LBFGS method can be competitive with its non-restarted 
counterpart in both noiseless and noisy settings. NLCG variants are less 
sensitive to the use of a restarting condition, but perform less favorably 
in a noisy setting.

Our analysis assumes that our gradient and function estimates satisfy 
deterministic accuracy conditions. Extending our approach to probabilistic 
accuracy properties, based on the notion of 
oracles~\cite{ASBerahas_LCao_KScheinberg_2021}, represents a promising avenue 
for future research. Investigating other line search techniques would also 
allow to hew closer to state-of-the-art LBFGS and NLCG techniques.

\paragraph{Acknowledgments}
We thank the coordinating editor and two anonymous reviewers for their 
feedback. In particular, we are grateful to one reviewer for pointing out the 
possibility of choosing $\alphapq$ independently of the Lipschitz constant.

Support for this research was partially provided by Agence Nationale de la 
Recherche through program ANR-19-P3IA-0001 (PRAIRIE 3IA Institute), by the FACE 
Foundation under the Thomas Jefferson Fund project ``Adaptive, Local and 
Innovative Algorithms for Stochastic Optimization (ALIAS)'', and by the 
Office of Naval Research under award N00014-24-1-2638.

\bibliographystyle{plain}
\bibliography{refsNCGrevised}

\begin{thebibliography}{10}

\bibitem{CAudet_WHare_2017}
C.~Audet and W.~Hare.
\newblock {\em Derivative-{F}ree and {B}lackbox {O}ptimization}.
\newblock Springer Series in Operations Research and Financial Engineering.
  Springer International Publishing, 2017.

\bibitem{ASBerahas_RHByrd_JNocedal_2019}
A.~S. {Berahas}, R.~H. {Byrd}, and J.~{Nocedal}.
\newblock Derivative-free optimization of noisy functions via quasi-newton
  methods.
\newblock {\em SIAM J. Optim.}, 29:965--993, 2019.

\bibitem{ASBerahas_LCao_KScheinberg_2021}
A.~S. {Berahas}, L.~Cao, and K.~Scheinberg.
\newblock Global convergence rate analysis of a generic line search algorithm
  with noise.
\newblock {\em SIAM J. Optim.}, 31:1489--1518, 2021.

\bibitem{EBerglund_JZhang_MJohansson_2025}
E.~Berglund, J.~Zhang, and M.~Johansson.
\newblock Soft quasi-{N}ewton: guaranteed positive definiteness by relaxing the
  secant constraint.
\newblock {\em Optim. Methods Softw.}, 40:783--812, 2025.

\bibitem{LBottou_FECurtis_JNocedal_2018}
L.~Bottou, F.~E. {Curtis}, and J.~Nocedal.
\newblock Optimization {M}ethods for {L}arge-{S}cale {M}achine {L}earning.
\newblock {\em SIAM Rev.}, 60:223--311, 2018.

\bibitem{byrd2012sample}
R.~H. Byrd, G.~M. Chin, J.~Nocedal, and Y.~Wu.
\newblock Sample size selection in optimization methods for machine learning.
\newblock {\em Math. Program.}, 134:127--155, 2012.

\bibitem{cao2024first}
L.~Cao, A.~S. Berahas, and K.~Scheinberg.
\newblock First-and second-order high probability complexity bounds for
  trust-region methods with noisy oracles.
\newblock {\em Math. Program.}, 207:55--106, 2024.

\bibitem{YCarmon_JCDuchi_OHinder_ASidford_2017a}
Y.~Carmon, J.~C. {Duchi}, O.~Hinder, and A.~Sidford.
\newblock ``{C}onvex until proven guilty": {D}imension-free acceleration of
  gradient descent on non-convex functions.
\newblock In {\em Proceedings of the International Conference on Machine
  Learning, August 2017, Sydney, Australia}, pages 654--663, 2017.

\bibitem{carter1991global}
R.~G. Carter.
\newblock On the global convergence of trust region algorithms using inexact
  gradient information.
\newblock {\em SIAM J. Numer. Anal.}, 28:251--265, 1991.

\bibitem{CCartis_PhRSampaio_PhLToint_2015}
C.~Cartis, Ph.~R. {Sampaio}, and Ph.~L. {Toint}.
\newblock Worst-case evaluation complexity of non-monotone gradient-related
  algorithms for unconstrained optimization.
\newblock {\em Optimization}, 64:1349--1361, 2015.

\bibitem{CCartis_KScheinberg_2018}
C.~Cartis and K.~Scheinberg.
\newblock Global convergence rate analysis of unconstrained optimization
  methods based on probabilistic models.
\newblock {\em Math. Program.}, 169:337--375, 2018.

\bibitem{RChanRenousLegoubin_CWRoyer_2022}
R.~{Chan--Renous-Legoubin} and C.~W. {Royer}.
\newblock A nonlinear conjugate gradient method with complexity guarantees and
  its application to nonconvex regression.
\newblock {\em Euro. J. Comput. Optim.}, 10:100044, 2022.

\bibitem{ARConn_KScheinberg_LNVicents_2009}
A.~R. {Conn}, K.~{Scheinberg}, and L.~N. {Vicente}.
\newblock {\em Introduction to derivative-free optimization}.
\newblock SIAM, 2009.

\bibitem{EDDolan_JJMore_2002}
E.~D. {Dolan} and J.~J. {Mor\'{e}}.
\newblock Benchmarking optimization software with performance profiles.
\newblock {\em Math. Program.}, 91:201--213, 2002.

\bibitem{JCDuchi_MIJordan_MJWainwright_AWibisono_2015}
J.~C. {Duchi}, M.~I. {Jordan}, M.~J. {Wainwright}, and A.~{Wibisono}.
\newblock Optimal rates for zero-order convex optimization: The power of two
  function evaluations.
\newblock {\em IEEE Transactions on Information Theory}, 61:2788--2806, 2015.

\bibitem{MFazel_RGe_SKakade_MMesbahi_2018}
M.~{Fazel}, R.~{Ge}, S.~{Kakade}, and M.~{Mesbahi}.
\newblock Global convergence of policy gradient methods for the linear
  quadratic regulator.
\newblock In {\em International Conference on Machine Learning}, pages
  1467--1476. PMLR, 2018.

\bibitem{NIMGould_DOrban_PhLToint_2015}
N.~I.~M. {Gould}, D.~Orban, and Ph.~L. {Toint}.
\newblock {CUTE}st: a constrained and unconstrained testing environment with
  safe threads.
\newblock {\em Comput. Optim. Appl.}, 60:545--557, 2015.

\bibitem{SGratton_PhLToint_2024}
S.~Gratton and Ph.~L. Toint.
\newblock {S2MPJ and CUTEst optimization problems for Matlab, Python and
  Julia}.
\newblock arXiv:2407.07812, 2024.

\bibitem{LGrippo_SLucidi_2005}
L.~Grippo and S.~Lucidi.
\newblock Convergence conditions, line search algorithms and trust region
  implementations for the {Polak-Ribi\`ere} conjugate gradient method.
\newblock {\em Optim. Methods Softw.}, 20:71--98, 2005.

\bibitem{WWHager_HZhang_2005}
W.~W. {Hager} and H.~Zhang.
\newblock A new conjugate gradient method with guaranteed descent and an
  efficient line search.
\newblock {\em SIAM J. Optim.}, 16:170--192, 2005.

\bibitem{WWHager_HZhang_2006a}
W.~W. {Hager} and H.~Zhang.
\newblock Algorithm 851: {CG\_DESCENT}, a conjugate gradient method with
  guaranteed descent.
\newblock {\em ACM Trans. Math. Software}, 32:113--137, 2006.

\bibitem{WWHager_HZhang_2006b}
W.~W. {Hager} and H.~Zhang.
\newblock A survey of nonlinear conjugate gradient methods.
\newblock {\em Pac. J. Optim.}, 2:35--58, 2006.

\bibitem{BIrwin_EHaber_2023}
B.~Irwin and E.~Haber.
\newblock Secant penalized {BFGS}: a noise robust quasi-{N}ewton method via
  penalizing the secant condition.
\newblock {\em Comput. Optim. Appl.}, 84:651--702, 2023.

\bibitem{BJin_KScheinberg_MXie_2024}
B.~Jin, K.~Scheinberg, and M.~Xie.
\newblock High probability complexity bounds for adaptive step search based on
  stochastic oracles.
\newblock {\em SIAM J. Optim.}, 34:2411--2439, 2024.

\bibitem{JLarson_MMenickelley_SMWild_2019}
J.~{Larson}, M.~{Menickelly}, and S.~M. {Wild}.
\newblock Derivative-free optimization methods.
\newblock {\em Acta Numerica}, 28:287--404, 2019.

\bibitem{JJMore_SMWild_2009}
J.J. {Mor\'{e}} and S.~M. {Wild}.
\newblock Benchmarking derivative-free optimization algorithms.
\newblock {\em SIAM J. Optim.}, 20:172--191, 2009.

\bibitem{JNocedal_SJWright_2006}
J.~Nocedal and S.~J. {Wright}.
\newblock {\em Numerical {O}ptimization}.
\newblock Springer Ser. Oper. Res. Financ. Eng. Springer-Verlag, New York,
  second edition, 2006.

\bibitem{CPaquette_KScheinberg_2020}
C.~Paquette and K.~Scheinberg.
\newblock A stochastic line search method with convergence rate analysis.
\newblock {\em SIAM J. Optim.}, 30:349--376, 2020.

\bibitem{RPasupathy_PGlynn_SGhosh_FSHashemi_2018}
R.~{Pasupathy}, P.~{Glynn}, S.~{Ghosh}, and F.~S. {Hashemi}.
\newblock On sampling rates in simulation-based recursions.
\newblock {\em SIAM J. Optim.}, 28(1):45--73, 2018.

\bibitem{SShashaani_FSHashemi_RPasupathy_2018}
S.~{Shashaani}, F.~S. {Hashemi}, and R.~{Pasupathy}.
\newblock Astro-df: A class of adaptive sampling trust-region algorithms for
  derivative-free stochastic optimization.
\newblock {\em SIAM J. Optim.}, 28(4):3145--3176, 2018.

\bibitem{shi2022noise}
H.-J.~M. {Shi}, Y.~Xie, R.~H. {Byrd}, and J.~Nocedal.
\newblock A noise-tolerant quasi-{N}ewton algorithm for unconstrained
  optimization.
\newblock {\em SIAM J. Optim.}, 32:29--55, 2022.

\bibitem{sun2023trust}
S.~Sun and J.~Nocedal.
\newblock A trust region method for noisy unconstrained optimization.
\newblock {\em Math. Program.}, 202:445--472, 2023.

\bibitem{SJWright_BRecht_2022}
S.~J. {Wright} and B.~Recht.
\newblock {\em Optimization for Data Analysis}.
\newblock Cambridge University Press, 2022.

\bibitem{xie2020analysis}
Y.~Xie, R.~H. {Byrd}, and J.~Nocedal.
\newblock Analysis of the {BFGS} method with errors.
\newblock {\em SIAM J. Optim.}, 30:182--209, 2020.

\end{thebibliography}

\begin{appendices}

\newpage
\section{Full restart statistics in the noisy setting}
\label{app:restarttables}

In Tables~\ref{tab:rstnoise1em8}-\ref{tab:rstnoise1em1}, bold numbers 
represent the pair $(p,\kappa)$ corresponding to the lowest percentage of 
restarts in the corresponding table, and italicized numbers represent the 
best overall values of $p$ (resp., $\kappad$) for different values of 
$\kappad$ (resp., $p$).

\begin{table}[ht]
{\footnotesize
\centering
\begin{tabular}{l|crrrrr}
\toprule
\multicolumn{7}{c}{NLCGr} \\
 \midrule
$p\backslash\kappad$ & &$10^{2}$ &$10^{3}$ &$10^{4}$ &$10^{5}$ &$\mathbf{10^{6}}$ \\
\cmidrule{1-7}
0~~~ & &59.04 &23.04 &11.28 &10.62 &\emph{10.58} \\
0.25 & &31.80 &11.10 &10.62 &10.60 &\emph{10.58} \\
0.50 & &11.27 &\textbf{10.40} &10.63 &10.61 &\emph{10.60} \\
\textbf{0.75} & &\emph{10.71} &\emph{10.48} &\emph{10.62} &\emph{10.60} &\emph{10.60} \\
1~~~ & &11.07 &10.58 &10.65 &10.61 &\emph{10.61} \\
\bottomrule
\end{tabular}
\quad 
\begin{tabular}{l|crrrrr}
\toprule
\multicolumn{7}{c}{LBFGSr} \\
 \midrule
$p\backslash\kappa$ & &$10^{2}$ &$10^{3}$ &$10^{4}$ &$10^{5}$ &$\mathbf{10^{6}}$ \\
 \cmidrule{1-7}
0~~~ & &46.18	&27.77	 &18.47	 &13.15	 &\emph{~6.13} \\
0.25 & &35.43	&22.17	 &15.58	 &~9.94	 &\emph{~5.18} \\
0.50 & &32.87	&19.07	 &12.43	 &~7.17	 &\emph{~3.81} \\
\textbf{0.75} & &\emph{34.94}	&\emph{16.99}	 &\emph{11.48}  &\emph{~4.62}  &\emph{\textbf{~3.54}} \\
1~~~ & &41.50	&17.26	 &~8.68	 &~4.89	 &\emph{~3.35} \\
\bottomrule
\end{tabular}}
\caption{Percentage of restarted iterations with $\epsf=10^{-8}$ (Left: NLCG; Right: LBFGS).}
\label{tab:rstnoise1em8}
\end{table}

\begin{table}[ht]
{\footnotesize
\centering
\begin{tabular}{l|crrrrr}
\toprule
\multicolumn{7}{c}{NLCGr} \\
 \midrule
$p\backslash\kappad$ & &$10^{2}$ &$10^{3}$ &$10^{4}$ &$\mathbf{10^{5}}$ &$10^{6}$ \\
 \cmidrule{1-7}
0~~~ & &16.30 &14.11 &13.81 &\emph{\textbf{13.72}} &\textbf{13.72} \\
0.25 & &14.82 &14.05 &13.78 &\emph{\textbf{13.72}} &13.77 \\
0.50 & &14.63 &13.99 &13.75 &\emph{13.73} &13.82 \\
\textbf{0.75} & &\emph{14.51} &\emph{13.95} &\emph{13.73} &\emph{13.77} &\emph{13.89} \\
1~~~ & &14.53 &13.90 &13.75 &\emph{13.86} &13.89 \\
\bottomrule
\end{tabular}
\quad 
\begin{tabular}{l|crrrrr}
\toprule
\multicolumn{7}{c}{LBFGSr} \\
 \midrule
$p\backslash\kappa$ & &$10^{2}$ &$10^{3}$ &$10^{4}$ &$10^{5}$ &$\mathbf{10^{6}}$ \\
\cmidrule{1-7}
0~~~ & &26.11	&16.63	 &~8.88	 &~5.83	 &\emph{~2.74} \\
0.25 & &23.10	&14.93	 &~8.35	 &~5.18	 &\textbf{~2.52} \\
0.50 & &21.31	&13.37	 &~8.07	 &~4.02	 &\emph{~2.53} \\
0.75 & &20.39	&12.32	 &~6.97  &~3.91  &\emph{~2.59} \\
\textbf{1~~~} & &\emph{20.07}	&\emph{11.91}	 &~\emph{5.83}	 &~\emph{3.98}	 &\emph{~2.62} \\
\bottomrule
\end{tabular}}
\caption{Percentage of restarted iterations with $\epsf=10^{-4}$ (Left: NLCG; Right: LBFGS).}
\label{tab:rstnoise1em4}
\end{table}

\begin{table}[ht]
{\footnotesize
\centering
\begin{tabular}{l|crrrrr}
\toprule
\multicolumn{7}{c}{NLCGr} \\
 \midrule
$p\backslash\kappad$ & &$10^{2}$ &$\mathbf{10^{3}}$ &$10^{4}$ &$10^{5}$ &$10^{6}$ \\
\cmidrule{1-7}
\textbf{0~~~} & &\emph{14.22} &\emph{13.82} &\emph{\textbf{13.71}} &\emph{13.82} &\emph{13.87} \\
0.25 & &14.16 &\emph{13.74} &13.75 &13.86 &13.94 \\
0.50 & &14.11 &\emph{13.72} &13.80 &13.88 &13.95 \\
0.75 & &14.05 &\emph{13.74} &13.84 &13.94 &13.94 \\
1~~~ & &14.01 &\emph{13.81} &13.88 &13.95 &13.94 \\
\bottomrule
\end{tabular}
\quad 
\begin{tabular}{l|crrrrr}
\toprule
\multicolumn{7}{c}{LBFGSr} \\
 \midrule
$p\backslash\kappa$ & &$10^{2}$ &$10^{3}$ &$10^{4}$ &$10^{5}$ &$\mathbf{10^{6}}$ \\
\cmidrule{1-7}
0~~~ & &17.62	&~7.56	 &~4.32	 &~1.75	 &\emph{~\textbf{1.18}} \\
0.25 & &16.38	&~8.22	 &~4.48	 &~1.92	 &\emph{~1.34} \\
\textbf{0.50} & &\emph{15.51}	&\emph{~8.37}	 &\emph{~4.43}	 &\emph{~2.15}	 &\emph{~1.57} \\
0.75 & &15.32	&~8.48	 &~4.27  &~2.43  &\emph{~1.71} \\
1~~~ & &14.81	&~8.78	 &~4.45	 &~2.53	 &\emph{~1.87} \\
\bottomrule
\end{tabular}}
\caption{Percentage of restarted iterations with $\epsf=10^{-2}$ (Left: NLCG; Right: LBFGS).}
\label{tab:rstnoise1em2}
\end{table}

\begin{table}[ht]
{\footnotesize
\centering
\begin{tabular}{l|crrrrr}
\toprule
\multicolumn{7}{c}{NLCGr} \\
 \midrule
$p\backslash\kappad$ & &$10^{2}$ &$\mathbf{10^{3}}$ &$10^{4}$ &$10^{5}$ &$10^{6}$ \\
 \cmidrule{1-7}
0~~~ & &16.22 &\emph{15.81} &15.83 &15.89 &15.93 \\
\textbf{0.25} & &\emph{16.08} &\emph{\textbf{15.80}} &\emph{15.88} &\emph{15.92} &\emph{15.92} \\
0.50 & &16.05 &\emph{15.83} &15.88 &15.93 &15.92 \\
0.75 & &16.09 &\emph{15.83} &15.91 &15.92 &15.92 \\
1~~~ & &16.08 &\emph{15.90} &15.93 &15.93 &15.92 \\
\bottomrule
\end{tabular}
\quad 
\begin{tabular}{l|crrrrr}
\toprule
\multicolumn{7}{c}{LBFGSr} \\
 \midrule
$p\backslash\kappa$ & &$10^{2}$ &$10^{3}$ &$10^{4}$ &$10^{5}$ &$\mathbf{10^{6}}$ \\
 \cmidrule{1-7}
\textbf{0~~~} & &\emph{~9.85}	&\emph{~5.02}	 &\emph{~3.21}	 &\emph{~1.29}	 &\emph{~\textbf{0.89}} \\
0.25 & &10.36	&~5.68	 &~3.25	 &~1.66	 &\emph{~1.07} \\
0.50 & &11.43	&~6.64	 &~3.49	 &~1.87	 &\emph{~1.39} \\
0.75 & &11.87	&~6.90	 &~3.78  &~2.24  &\emph{~1.50} \\
1~~~ & &12.13	&~7.09	 &~4.20	 &~2.46	 &\emph{~1.53} \\
 \bottomrule
\end{tabular}}
\caption{Percentage of restarted iterations with $\epsf=10^{-1}$ (Left: NLCG; Right: LBFGS).}
\label{tab:rstnoise1em1}
\end{table}

\end{appendices}

\end{document}